\documentclass[11pt, DIV10,a4paper]{article}
\usepackage{color}
\usepackage{float}
\usepackage[bf,small]{caption2}
\usepackage{comment}
\usepackage{amsmath}
\usepackage{bigints}
\usepackage{rotating, booktabs}
\usepackage{amsopn}
\usepackage{amsfonts}
\usepackage{amsthm}
\usepackage{amssymb}
\usepackage{amsbsy}
\usepackage{multirow}
\usepackage{slashbox}
\usepackage{natbib}
\usepackage{tocbibind}
\usepackage[a4paper,colorlinks,breaklinks,bookmarksopen,bookmarksnumbered]{hyperref}
\usepackage[refpage]{nomencl}
\usepackage{makeidx}
\usepackage{pst-all}
\usepackage{epsfig,psfrag}
\usepackage{graphicx}
\usepackage{amsmath}
\usepackage{epstopdf}
\usepackage{tabularx}
\usepackage{subfigure}
\usepackage{setspace}
\usepackage[mathscr]{euscript}
\usepackage[margin=1in]{geometry}
\usepackage{xr-hyper}
\usepackage{amsfonts} 
\usepackage{calc}

\newcommand{\ueq}[1][]{%
  \if\relax\detokenize{#1}\relax
    \sbox0{$\underbrace{=}_{}$}%
    \mathrel{\mathmakebox[\wd0]{=}}
  \else
    \mathrel{\underbrace{=}_{\mathclap{#1}}}
  \fi}

\newcommand{\bzero}{\boldsymbol{0}}

\newcommand {\ctn}{\cite}

\usepackage{url}
\newcommand{\e}{\ensuremath{\epsilon}}

\newcommand {\btheta}{\mbox{$\theta$}}

\newcommand {\bdelta}{\mbox{\boldmath $\delta$}}

\newcommand {\bTheta}{\mbox{$\Theta$}}
\newcommand {\boeta}{\mbox{\boldmath $\eta$}}

\newcommand{\bI}{\mathbf I}

\newcommand{\bx}{\mathbf x}
\newcommand{\bX}{\mathbf X}
\newcommand{\by}{\mathbf y}
\newcommand{\bY}{\mathbf Y}
\newcommand{\bz}{\mathbf z}

\newtheorem{theorem}{Theorem}

\numberwithin{equation}{section}
\numberwithin{algo}{section}
\numberwithin{table}{section}
\numberwithin{figure}{section}

\newtheorem{lemma}{Lemma}[section]

\normalsize

\begin{document}

\title{\textbf{Posterior Convergence of Gaussian and General Stochastic Process Regression Under Possible Misspecifications}}
\author{Debashis Chatterjee$^{\dag}$ and Sourabh Bhattacharya$^{\dag, +}$ }
\date{}
\maketitle
\begin{center}
$^{\ddag}$ Indian Statistical Institute\\
$+$ Corresponding author:  \href{mailto: bhsourabh@gmail.com}{bhsourabh@gmail.com}
\end{center}

\begin{abstract}
In this article, we investigate posterior convergence in nonparametric regression models where the unknown regression function is modeled by some
appropriate stochastic process. In this regard, we consider two setups. The first setup is based on Gaussian processes, 
where the covariates are either random or non-random and the noise may be either normally or double-exponentially distributed. 
In the second setup, we assume that the underlying regression function is modeled by some reasonably smooth, but unspecified stochastic process satisfying
reasonable conditions. The distribution of the noise is also left unspecified, but assumed to be thick-tailed.
As in the previous studies regarding the same problems, we do not assume that the truth lies in the postulated
parameter space, thus explicitly allowing the possibilities of misspecification. We exploit the general results of \ctn{Shalizi09} for our purpose and 
establish not only posterior consistency, but also the rates at which the posterior probabilities converge, which turns out to be the Kullback-Leibler divergence rate. 
We also investigate the more familiar posterior convergence rates. Interestingly, we show that the posterior predictive distribution can accurately 
approximate the best possible predictive distribution in the sense that the Hellinger distance, as well as the total variation distance between the two 
distributions can tend to zero, in spite of misspecifications.
\\[2mm]
{\bf Keywords:} {\it Double exponential distribution; Gaussian process; Infinite dimension; Kullback-Leibler divergence rate; Misspecification; Posterior convergence.}
\end{abstract}

\section{Introduction}
\label{sec:intro}

In statistics, either frequentist or Bayesian, nonparametric regression plays a very significant role. The frequentist nonparametric literature, however, is substantially
larger than the Bayesian counterpart. Here we cite the books \ctn{Schimek13}, \ctn{Hardle12}, \ctn{Efro08}, \ctn{Tak06}, \ctn{Wu06}, \ctn{Eubank99}, \ctn{Green93} and \ctn{Hardle90}, 
among a large number of books on frequentist nonparametric regression. 
The Bayesian nonparametric literature, which is relatively young but flourishing in the recent times (see, for example, \ctn{Ghosal17}, \ctn{Muller15}, \ctn{Dey12},
\ctn{Hjort10}, \ctn{Ghosh03}), offers much broader scope for
interesting and innovative research. 

The importance of Gaussian processes in nonparametric statistical modeling, particularly in the Bayesian context, is undeniable. It is widely used in density estimation 
(\ctn{Lenk88}, \ctn{Lenk91}, \ctn{Lenk03}),
nonparametric regression (\ctn{Ras06}), spatial data modeling (\ctn{Cressie93}, \ctn{Banerjee14}), machine learning (\ctn{Ras06}), emulation of computer models (\ctn{Santner03}), 
to name a few areas. 
Although applications of Gaussian processes have 
received and continue to receive much attention, 
in the recent years there seems to be a growing interest among researchers in the theoretical properties of approaches based on Gaussian processes.
Specifically, investigation of posterior convergence of Gaussian process based approaches has turned out to be an important undertaking. In this respect, contributions are made by
\ctn{Choi07}, \ctn{Vaart08}, \ctn{Vaart09}, \ctn{Vaart11}, \ctn{Knapik11}, \ctn{Vollmer13}, \ctn{yang17}, \ctn{Knapik18}. \ctn{Choi07} address posterior consistency in Gaussian process
regression, while the others also attempt to provide the rates of posterior convergence. However, the rates are so far computed under the assumption that the error distribution is normal
and the error variance is either known, or if unknown, can be given a prior, but on a compact support bounded away from zero. 

General priors for the regression function or thick-tailed noise distributions seemed to have received
less attention. The asymptotic theory for such frameworks is even rare, \ctn{Choi09} being an important exception. 
As much as we are aware of, rates of convergence are not available for nonparametric 
regression with general stochastic process prior on the regression function and thick-tailed noise distributions. Another important issue which seems to have received 
less attention in the literature, is the case of misspecified models.
We are not aware of any published asymptotic theory pertaining to misspecifications in nonparametric regression, for either Gaussian or non-Gaussian processes with either normal or non-normal errors.

In this article, we consider both Gaussian and general stochastic process regression under the same setups as \ctn{Choi07} and \ctn{Choi09}, respectively, 
assuming that the covariates may be either random or non-random. For the Gaussian process setup we consider both
normal and double-exponential distribution for the error, with unknown error variance. In the general context, we assume non-Gaussian noise with unknown scale 
parameter supported on the entire positive part of the real line.
Based on the general theory of posterior convergence provided in \ctn{Shalizi09},
we establish posterior convergence theories for both the setups. We allow the case of misspecified models, that is, if the true regression function and the true error variance
are not even supported by the prior. Our approach also enables us to show that the relevant posterior probabilities converge at the Kullback-Leibler (KL) divergence rate, and that
the posterior convergence rate with respect to the KL-divergence is just slower than $n^{-1}$, $n$ being the number of observations. We further show that even in the case of misspecification,
the posterior predictive distribution can approximate the best possible predictive distribution adequately, in the sense that
the Hellinger distance, as well as the total variation distance between the two distributions can tend to zero.
In Section \ref{subsec:shalizi_briefing} we provide a brief overview and intuitive explanation of the main assumptions and results of Shalizi, which we exploit in this article.
The details are provided in Section \ref{sec:shalizi} of the supplement. The results of Shalizi are based on seven assumptions, which we refer to as (S1) -- (S7) in this article.

\subsection{A briefing of the main results of Shalizi}
\label{subsec:shalizi_briefing}
Let $\bY_n=(Y_1,\ldots,Y_n)^T$, and let $f_{\theta}(\bY_n)$ and $f_{\theta_0}(\bY_n)$ denote the observed and the true likelihoods respectively, under the given value of the parameter $\theta$
and the true parameter $\theta_0$. We assume that $\theta\in\Theta$, where $\Theta$ is the (often infinite-dimensional) parameter space. However, we {\it do not} 
assume that $\theta_0\in\Theta$, thus allowing misspecification.
The key ingredient associated with Shalizi's approach to proving convergence of the posterior distribution of $\theta$ is to show that the asymptotic equipartition property holds.
To elucidate, let us consider the following likelihood ratio:
\begin{equation*}
R_n(\theta)=\frac{f_{\theta}(\bY_n)}{f_{\theta_0}(\bY_n)}.
\end{equation*}
Then, to say that for each $\theta\in\Theta$, the generalized or relative asymptotic equipartition property holds, we mean
\begin{equation}
\underset{n\rightarrow\infty}{\lim}~\frac{1}{n}\log R_n(\theta)=-h(\theta),
\label{eq:equipartition}
\end{equation}
almost surely, where
$h(\theta)$ is the KL-divergence rate given by
\begin{equation*}
h(\theta)=\underset{n\rightarrow\infty}{\lim}~\frac{1}{n}E_{\theta_0}\left(\log\frac{f_{\theta_0}(\bY_n)}{f_{\theta}(\bY_n)}\right),
\end{equation*}
provided that it exists (possibly being infinite), where $E_{\theta_0}$ denotes expectation with respect to the true model.
Let
\begin{align}
h\left(A\right)&=\underset{\theta\in A}{\mbox{ess~inf}}~h(\theta);\notag\\ 
J(\theta)&=h(\theta)-h(\Theta);\notag\\ 
J(A)&=\underset{\theta\in A}{\mbox{ess~inf}}~J(\theta).\notag 
\end{align}
Thus, $h(A)$ can be roughly interpreted as the minimum KL-divergence between the postulated and the true model over the set $A$. If $h(\Theta)>0$, this indicates
model misspecification. However, as we shall show, model misspecification need not always imply that $h(\Theta)>0$. For $A\subset\Theta$, $h(A)>h(\Theta)$, so that
$J(A)>0$.

As regards the prior, it is required to construct an appropriate sequence of sieves $\mathcal G_n$ such that $\mathcal G_n\rightarrow\Theta$ and $\pi(\mathcal G^c_n)\leq\alpha\exp(-\beta n)$,
for some $\alpha>0$. 

With the above notions, verification of (\ref{eq:equipartition}) along with several other technical conditions ensure that for any $A\subseteq\Theta$ such that $\pi(A)>0$, 
\begin{equation}
\underset{n\rightarrow\infty}{\lim}~\pi(A|\bY_n)=0,
\label{eq:post_conv1}
\end{equation}
almost surely, provided that $h(A)>h(\Theta)$. Under mild assumptions, it also holds that
\begin{equation}
\underset{n\rightarrow\infty}{\lim}~\frac{1}{n}\log\pi(A|\bY_n)=-J(A),
\label{eq:post_conv2}
\end{equation}
almost surely, where $\pi(\cdot|\bY_n)$ denotes the posterior distribution of $\theta$ given $\bY_n$.
With respect to (\ref{eq:post_conv1}) note that $h(A)>h(\Theta)$ implies positive KL-divergence in $A$, even if $h(\Theta)=0$. In other words, $A$ is the set in which the postulated model
fails to capture the true model in terms of the KL-divergence. Hence, expectedly, the posterior probability of that set converges to zero. The result (\ref{eq:post_conv2})
asserts that the rate at which the posterior probability of $A$ converges to zero is about $\exp(-nJ(A))$. From the above results it is clear that the posterior concentrates
on sets of the form $N_{\epsilon}=\left\{\theta:h(\theta)\leq h(\Theta)+\epsilon\right\}$, for any $\epsilon>0$.

As regards the rate of posterior convergence, let $N_{\epsilon_n}=\left\{\theta:h(\theta)\leq h(\Theta)+\epsilon_n\right\}$, where $\epsilon_n\rightarrow 0$ such that $n\epsilon_n\rightarrow\infty$.
Then under an additional technical assumption it holds, almost surely, that
\begin{equation}
\underset{n\rightarrow\infty}{\lim}~\pi\left(N_{\epsilon_n}|\bY_n\right)=1.
\label{eq:conv_rate}
\end{equation}

Moreover, it was shown by Shalizi that the squares of the Hellinger and the total variation distances between the posterior predictive distribution 
and the best possible predictive distribution under the truth, are asymptotically almost surely bounded above by $h(\Theta)$ and $4h(\Theta)$, respectively. In other words, if $h(\Theta)=0$,
then this entails very accurate approximation of the true predictive distribution by the posterior predictive distribution.

The rest of our article is structured as follows. We treat the Gaussian process regression with normal and double exponential errors in Section \ref{sec:gp_setup}.
Specifically, our assumptions regarding the model and discussion of the assumptions are presented in Section \ref{subsec:assumptions}. 
In Section \ref{subsec:postconv_normal_de_errors} we present our main results of posterior convergence, 
along with the summary of the verification of Shalizi's assumptions, for the Gaussian process setup. 
The complete details are provided in Sections \ref{sec:verification} and \ref{sec:de} of the supplement. 
We deal with rate of convergence and model misspecification issue for Gaussian process regression in Sections \ref{subsec:rate} and \ref{subsec:miss}, respectively. 

The case of general stochastic process regression with thick tailed error distribution is taken up in Section \ref{sec:nongp_setup}. The assumptions with their discussion
are provided in Section \ref{subsec:nongp_assumptions}, the main posterior results are presented in Section \ref{subsec:nongp_postconv_normal_de_errors}, and
Section \ref{subsec:nongp_rate} addresses the rate of convergence and model misspecification issue.
Finally, we make concluding remarks in Section \ref{sec:conclusion}.
The relevant details are provided in Section \ref{sec:nongp_verification} of the supplement.

\section{The Gaussian process regression setup}
\label{sec:gp_setup}
As in \ctn{Choi07}, we consider the following model:
\begin{align}
y_i&=\eta(\bx_i)+\epsilon_i;~i=1,\ldots,n;\label{eq:model1}\\
\eta(\cdot) &\sim GP\left(\mu(\cdot),c(\cdot,\cdot)\right);\label{eq:gp1}\\
\sigma &\sim\pi_{\sigma}(\cdot).
\label{eq:prior_sigma}
\end{align}
In (\ref{eq:gp1}), $GP\left(\mu(\cdot),c(\cdot,\cdot)\right)$ stands for Gaussian process with mean function $\mu(\cdot)$ and positive definite covariance function 
$cov(\eta(\bx_1),\eta(\bx_2))=c(\bx_1,\bx_2)$, for any $\bx_1,\bx_2\in\mathcal X$, where $\mathcal X$ is the domain of $\eta$. 

As in \ctn{Choi07} we assume two separate distributions for the errors $\epsilon_i$, independent zero-mean normal with variance $\sigma^2$ which we denote by $N(0,\sigma^2)$ and
independent double exponential distribution with median $0$ and scale parameter $\sigma$ with density
\begin{equation*}
f(\epsilon)=\frac{1}{2\sigma}\exp\left(-\frac{|\epsilon|}{\sigma}\right);~\epsilon\in\mathbb R.
\end{equation*}
We denote the double exponential distribution by $DE(0,\sigma)$.

In our case, let $\theta=(\eta,\sigma)$ be the infinite-dimensional parameter associated with our Gaussian process model and let $\theta_0=(\eta_0,\sigma_0)$ be the true
(infinite-dimensional) parameter. Let $\Theta$ denote the infinite-dimensional parameter space. 
%

\subsection{Assumptions and their discussions}
\label{subsec:assumptions}
Regarding the model and the prior, we make the following assumptions:
\begin{itemize}
\item[(A1)]$\mathcal X$ is a compact, $d$-dimensional space, for some finite $d\geq 1$, equipped with a suitable metric.
\item[(A2)] The functions $\eta$ are continuous on $\mathcal X$ and for such functions
the limit 
\begin{equation}
\eta'_j(\bx)=\frac{\partial \eta(\bx)}{\partial x_j}=\underset{h\rightarrow 0}{\lim}~\frac{\eta(\bx+h\bdelta_j)-\eta(\bx)}{h}
\label{eq:partial_derivative1}
\end{equation}
exists for each $\bx\in\mathcal X$, and is continuous on $\mathcal X$, for $j=1,\ldots,d$. 
In the above, $\bdelta_j$ is the $d$-dimensional vector where the $j$-th element is 1 and all the other elements are zero.
We denote the above class of functions by $\mathcal C'(\mathcal X)$.

\item[(A3)] We assume the following for the covariates $\bx_i$, accordingly as they are considered an observed random sample, or non-random.

\begin{enumerate}
\item[(i)]
$\left\{\bx_i:i=1,2,\ldots\right\}$ is an observed sample associated with an $iid$ sequence associated with some probability measure $Q$, supported on $\mathcal X$, which is independent
of $\left\{\e_i:i=1,2,\ldots\right\}$. 

\item[(ii)] $\left\{\bx_i:i=1,2,\ldots\right\}$ is an observed non-random sample. In this case, we consider 
a specific partition of the $d$-dimensional space $\mathcal X$ into $n$ subsets such that
each subset of the partition contains at least one $\bx\in\left\{\bx_i:i=1,2,\ldots\right\}$ and has Lebesgue measure $L/n$, for some $L>0$.  

\end{enumerate}

\item[(A4)] Regarding the prior for $\sigma$, we assume that for large enough $n$,
	$$\pi_{\sigma}\left(\exp(-\left(\beta n\right)^{1/4})\leq\sigma\leq\exp(\left(\beta n\right)^{1/4})\right)\geq 1-c_{\sigma}\exp(-\beta n),$$ for $c_{\sigma}>0$ and $\beta>2h\left(\Theta\right)$.

\item[(A5)] The true regression function $\eta_0$ satisfies $\|\eta_0\|\leq\kappa_0<\infty$. We {\it do not} assume that $\eta_0\in\mathcal C'(\mathcal X)$. For random covariate
$\bX$, we assume that $\eta_0(\bX)$ is measurable.



\end{itemize}

\subsubsection{Discussion of the assumptions}
\label{subsubsec:ass_diss}
The compactness assumption on $\mathcal X$ in Assumption (A1) guarantees that continuous functions on $\mathcal X$ have finite sup-norms. Here, by sup-norm of any function $f$ on $\mathcal X$,
we mean $\|f\|=\underset{\bx\in\mathcal X}{\sup}~|f(\bx)|$. Hence, our Gaussian process prior on $\eta$,
which gives probability one to continuously differentiable functions, also ensures that $\|\eta\|<\infty$, almost surely. Compact support of the functions is commonplace in the Gaussian
process literature; see, for example, \ctn{Cramer67}, \ctn{Adler81}, \ctn{Adler07}, \ctn{Choi07}.
The metric on $\mathcal X$ is necessary for partitioning $\mathcal X$ in the case of non-random covariates. 

Condition (A2) is required for constructing appropriate sieves for proving our posterior convergence results. In particular, this is required to ensure that $\eta$ is Lipschitz continuous
in the sieves. Since a function is Lipschitz if and only if its partial derivatives are bounded, this serves our purpose, as continuity of the partial derivatives of $\eta$ guarantees
boundedness in the compact domain $\mathcal X$.
Conditions guaranteeing the above continuity and smoothness properties required by (A2) must also be reflected in the underlying Gaussian process prior for $\eta$. 
The relevant conditions can be found in \ctn{Cramer67}, \ctn{Adler81} and \ctn{Adler07}, which we assume in our case. 
In particular, these require adequate smoothness assumptions on the mean function $\mu(\cdot)$ and the covariance function $c(\cdot,\cdot)$ of the Gaussian process prior. 
It follows that $\eta'_j$; $j=1,\ldots,d$, are also Gaussian processes. 
It clearly holds that $\mu(\cdot)$ and its partial derivatives also have finite sup-norms. 

As regards (A3) (i), thanks to the strong law of large numbers (SLLN), given any $\eta$ in the complement of some null set with respect to the prior, and given any sequence 
$\left\{\bx_i:i=1,2,\ldots\right\}$ 
this assumption ensures that for any $\nu>0$, as $n\rightarrow\infty$,
\begin{equation}
\frac{1}{n}\sum_{i=1}^n\left|\eta(\bx_i)-\eta_0(\bx_i)\right|^{\nu}\rightarrow\int_{\mathcal X}\left|\eta(\bx)-\eta_0(\bx)\right|^{\nu}dQ(\bX)
=E_{\bX}\left|\eta(\bX)-\eta_0(\bX)\right|^{\nu}~\mbox{(say)},
\label{eq:a3_1}
\end{equation}
where $Q$ is some probability measure supported on $\mathcal X$.

Condition (A3) (ii) ensures that $\frac{1}{n}\sum_{i=1}^n\left|\eta(\bx_i)-\eta_0(\bx_i)\right|^{\nu}$ is a particular Riemann sum and hence (\ref{eq:a3_1}) holds with 
$Q$ being the Lebesgue measure on $\mathcal X$. 
We continue to denote the limit in this case by $E_{\bX}\left[\eta(\bX)-\eta_0(\bX)\right]^{\nu}$. 

In the light of (\ref{eq:a3_1}), condition (A3) will play important role in establishing the equipartition property, for both Gaussian and double exponential errors. 
Another important role of this condition is to ensure consistency of the posterior predictive distribution, in spite of some misspecifications.

Condition (A4) ensures that the prior probabilities of the complements of the sieves are exponentially small. Such a requirement is common to most Bayesian asymptotic theories.

The essence of (A5) is to allow misspecification of the prior for $\eta$ in a way that the true regression function is not even supported by the prior, even though it has finite sup-norm.
In contrast, \ctn{Choi07} assumed that $\eta_0$ has continuous first-order partial derivatives. 
The assumption of measurability of $\eta_0(\bX)$ is a very mild technical condition.

Let $\Theta=\mathcal C'(\mathcal X)\times\mathbb R^+$ denote the infinite-dimensional parameter space for our Gaussian process model.

\subsection{Posterior convergence of Gaussian process regression under normal and double exponential errors}
\label{subsec:postconv_normal_de_errors}

In this section we provide a summary of our results leading to posterior convergence of Gaussian process regression when the errors are assumed to be either normal or double exponential. The details
are provided in the supplement. The key results associated with the asymptotic equipartition property are provided in Lemma \ref{lemma:lemma1} and Theorem \ref{theorem:theorem1}, the proofs
of which are provided in the supplement in the context of detailed verification of Shalizi's assumptions.
\begin{lemma}
\label{lemma:lemma1}
Under the Gaussian process model and conditions (A1) and (A3), the KL-divergence rate $h(\theta)$ exists for $\theta\in\Theta$, and is given by
\begin{equation}
h(\theta)=\log\left(\frac{\sigma}{\sigma_0}\right)-\frac{1}{2}+\frac{\sigma^2_0}{2\sigma^2}+\frac{1}{2\sigma^2}E_\bX\left[\eta(\bX)-\eta_0(\bX)\right]^2,
\label{eq:h}
\end{equation}
for the normal errors, and
\begin{equation}
h(\theta)=\log\left(\frac{\sigma}{\sigma_0}\right)-1+\frac{1}{\sigma}E_{\bX}\left|\eta(\bX)-\eta_0(\bX)\right|
+\frac{\sigma_0}{\sigma}E_{\bX}\left[\exp\left(-\frac{|\eta(\bX)-\eta_0(\bX)|}{\sigma_0}\right)\right],
\label{eq:h_de}
\end{equation}
for the double exponential errors.
\end{lemma}

\begin{theorem}
\label{theorem:theorem1}
Under the Gaussian process model with normal and double exponential errors and conditions (A1) and (A3), the asymptotic equipartition property holds, and is given by
\begin{equation*}
\underset{n\rightarrow\infty}{\lim}~\frac{1}{n}\log R_n(\theta)=-h(\theta),~\mbox{almost surely}.
\end{equation*}
The convergence is uniform on any compact subset of $\Theta$.
\end{theorem}

Lemma \ref{lemma:lemma1} and Theorem \ref{theorem:theorem1} ensure that conditions (S1) -- (S3) of Shalizi hold, and (S4) holds since $h(\theta)$ is almost surely
finite. We construct the sieves $\mathcal G_n$ as  
\begin{equation}
	\mathcal G_n=\left\{\left(\eta,\sigma\right):\|\eta\|\leq\exp(\left(\beta n\right)^{1/4}),\exp(-\left(\beta n\right)^{1/4})\leq\sigma\leq\exp(\left(\beta n\right))^{1/4}),
	\|\eta'_j\|\leq\exp(\left(\beta n\right)^{1/4});j=1,\ldots,d\right\}. 
\label{eq:G}
\end{equation}
It follows that $\mathcal G_n\rightarrow\Theta$ as $n\rightarrow\infty$ and the properties of the Gaussian processes $\eta$, $\eta'$,
together with (A4) ensure that $\pi(\mathcal G^c_n)\leq\alpha\exp(-\beta n)$, for some $\alpha>0$.
This result, continuity of $h(\theta)$, compactness of $\mathcal G_n$ and the uniform convergence result of Theorem \ref{theorem:theorem1}, together ensure (S5).

Now observe that the aim of assumption (S6) is to ensure that (see the proof of Lemma 7 of \ctn{Shalizi09}) 
for every $\varepsilon>0$ and for all $n$ sufficiently large,
\begin{equation*}
\frac{1}{n}\log\int_{\mathcal G_n}R_n(\btheta)d\pi(\btheta)\leq -h\left(\mathcal G_n\right)+\varepsilon,~\mbox{almost surely}.
\end{equation*}
Since $h\left(\mathcal G_n\right)\rightarrow h\left(\bTheta\right)$ as $n\rightarrow\infty$, it is enough to verify that for every $\varepsilon>0$ and for all $n$ sufficiently large,
\begin{equation}
\frac{1}{n}\log\int_{\mathcal G_n}R_n(\btheta)d\pi(\btheta)\leq -h\left(\bTheta\right)+\varepsilon,~\mbox{almost surely}.
\label{eq:s6_1}
\end{equation}
In this regard, first observe that 
\begin{align}
\frac{1}{n}\log\int_{\mathcal G_n}R_n(\btheta)d\pi(\btheta)&\leq\frac{1}{n}\log\left[\underset{\btheta\in\mathcal G_n}{\sup}~R_n(\btheta)\pi(\mathcal G_n)\right]\notag\\
&=\frac{1}{n}\log\left[\underset{\btheta\in\mathcal G_n}{\sup}~R_n(\btheta)\right]+\frac{1}{n}\log\pi(\mathcal G_n)\notag\\
&=\underset{\btheta\in\mathcal G_n}{\sup}~\frac{1}{n}\log R_n(\btheta)+\frac{1}{n}\log\pi(\mathcal G_n)\notag\\
&\leq \frac{1}{n}\underset{\btheta\in\mathcal G_n}{\sup}~\log R_n(\btheta),
\label{eq:s6_2}
\end{align}
where the last inequality holds since $\frac{1}{n}\log\pi(\mathcal G_n)\leq 0$. Now, letting $\mathcal S=\left\{\btheta:h(\btheta)\leq\kappa\right\}$, where
$\kappa>h\left(\bTheta\right)$ is large as desired,
\begin{align}
\underset{\btheta\in\mathcal G_n}{\sup}~\frac{1}{n}\log R_n(\btheta)
&\leq\underset{\btheta\in\bTheta}{\sup}~\frac{1}{n}\log R_n(\btheta)
=\underset{\btheta\in\mathcal S\cup\mathcal S^c}{\sup}~\frac{1}{n}\log R_n(\btheta)
\notag\\
&\leq \max\left\{\underset{\btheta\in\mathcal S}{\sup}~\frac{1}{n}\log R_n(\btheta),\underset{\btheta\in\mathcal S^c}{\sup}~\frac{1}{n}\log R_n(\btheta)\right\}.
\label{eq:s6_3}
\end{align}
In Sections \ref{subsubsec:S5_3} and \ref{subsec:S5_de} we have proved continuity of $h(\theta)$ for Gaussian and double exponential errors, respectively. 
Now observe that $\|\eta\|\leq\|\eta-\eta_0\|+\|\eta_0\|$, so that $\|\eta\|\rightarrow\infty$ implies $\|\eta-\eta_0\|\rightarrow\infty$ (since $\|\eta_0\|<\infty$).
Hence, for each $\eta$, there exists a subset $\mathcal X_{\eta}$ of $\mathcal X$ depending upon $\eta$ such that $Q\left(\mathcal X_{\eta}\right)>0$ and 
$\underset{x\in\mathcal X_{\eta}}{\sup}~|\eta(x)-\eta_0(x)|\rightarrow\infty$ as $\|\eta\|\rightarrow\infty$. It then follows that $E_{\bX}\left|\eta(\bX)-\eta_0(\bX)\right|\rightarrow\infty$ 
and $E_{\bX}\left(\eta(\bX)-\eta_0(\bX)\right)^2\rightarrow\infty$ as $\|\eta\|\rightarrow\infty$. Hence observe that $\|\theta\|\rightarrow\infty$ if
$\sigma\rightarrow \infty$ and $\|\eta\|\rightarrow\infty$, or if $\sigma$ tends to zero or some non-negative constant and $\|\eta\|\rightarrow\infty$. In both the cases
$h(\btheta)\rightarrow\infty$, for both Gaussian and double exponential errors. 
In other words,
$h(\btheta)$ is a continuous coercive function. Hence, $\mathcal S$ is a compact set (see, for example, \ctn{Lange10}). 
Now note that for any two real valued functions $f$ and $g$, and for any set $A$ in the domain of $f$ and $g$, 
$\underset{x\in A}{\sup}~f(x)=\underset{x\in A}{\sup}~\left[(f(x)-g(x))+g(x)\right]\leq\underset{x\in A}{\sup}~(f(x)-g(x))+\underset{x\in A}{\sup}~g(x)$,
so that 
$\underset{x\in A}{\sup}~f(x)-\underset{x\in A}{\sup}~g(x)\leq \underset{x\in A}{\sup}~(f(x)-g(x))\leq \underset{x\in A}{\sup}~\left|f(x)-g(x)\right|$.
Interchanging the roles of $f$ and $g$, we obtain $-\left[\underset{x\in A}{\sup}~f(x)-\underset{x\in A}{\sup}~g(x)\right]\leq \underset{x\in A}{\sup}~\left|f(x)-g(x)\right|$,
so that combining the above two inequalities lead to 
$\left|\underset{x\in A}{\sup}~f(x)-\underset{x\in A}{\sup}~g(x)\right|\leq \underset{x\in A}{\sup}~\left|f(x)-g(x)\right|$.
In our case, replacing $A$, $f(x)$ and $g(x)$ with $\mathcal S$, $\frac{1}{n}\log R_n(\btheta)$ and $-h(\btheta)$, respectively, we obtain using uniform convergence
of Theorem \ref{theorem:theorem1}, that
\begin{equation}
\underset{\btheta\in\mathcal S}{\sup}~\frac{1}{n}\log R_n(\btheta)\rightarrow \underset{\btheta\in\mathcal S}{\sup}~-h(\btheta)=-h\left(\mathcal S\right),~\mbox{almost surely, as}~n\rightarrow\infty.
\label{eq:s6_4}
\end{equation}

We now show that 
\begin{equation}
\underset{\btheta\in\mathcal S^c}{\sup}~\frac{1}{n}\log R_n(\btheta)\leq -h\left(\bTheta\right)~\mbox{almost surely, as}~ n\rightarrow\infty.
\label{eq:s6_5}
\end{equation}
First note that
if $\underset{\btheta\in\mathcal S^c}{\sup}~\frac{1}{n}\log R_n(\btheta)> -h\left(\bTheta\right)$ infinitely often, then $\frac{1}{n}\log R_n(\btheta)> -h\left(\bTheta\right)$
for some $\btheta\in\mathcal S^c$ infinitely often. But $\frac{1}{n}\log R_n(\btheta)> -h\left(\bTheta\right)$ if and only if 
$
\frac{1}{n}\log R_n(\btheta)+h(\btheta)> h(\btheta)-h\left(\bTheta\right),~\mbox{for}~\btheta\in\mathcal S^c. 
$
Hence, if we can show that  
\begin{equation}
P\left(\left|\frac{1}{n}\log R_n(\btheta)+h(\btheta)\right|> \kappa-h\left(\bTheta\right),~\mbox{for}~\btheta\in\mathcal S^c~\mbox{infinitely often}\right)=0, 
\label{eq:s6_6}
\end{equation}
then (\ref{eq:s6_5}) will be proved. We use the Borel-Cantelli lemma to prove (\ref{eq:s6_6}). In other words, we prove in the supplement, in the context of verifying condition (S6) of Shalizi, that
\begin{theorem}
\label{theorem:theorem3}
For both normal and double exponential errors, under (A1)--(A5), it holds that 
\begin{equation}
\sum_{n=1}^{\infty}\int_{\mathcal S^c}P\left(\left|\frac{1}{n}\log R_n(\btheta)+h(\btheta)\right|> \kappa-h\left(\bTheta\right)\right)d\pi(\btheta)<\infty. 
\label{eq:s6_7}
\end{equation}
\end{theorem}
%
Since $h(\btheta)$ is continuous, (S7) holds trivially. 
In other words, all the assumptions (S1)--(S7) are satisfied for Gaussian process regression, for both normal and double exponential errors. 
Formally, our results lead to the following theorem.
\begin{theorem}
\label{theorem:gp2}
Assume the Gaussian process regression model where the errors are either normally or double-exponentially distributed. 
Then under the conditions (A1) -- (A5), (\ref{eq:post_conv1}) holds. 
Also, for any measurable set $A$ with $\pi(A)>0$, if $\beta>2h(A)$, where $h$ is given by (\ref{eq:h}) for normal errors and (\ref{eq:h_de}) for double-exponential errors, 
or if $A\subset\cap_{k=n}^{\infty}\mathcal G_k$ for some $n$, 
where $\mathcal G_k$ is given by (\ref{eq:G}), then (\ref{eq:post_conv1}) and
(\ref{eq:post_conv2}) hold.
\end{theorem}

\subsection{Rate of convergence}
\label{subsec:rate}
Shalizi considered the set $N_{\epsilon_n}=\left\{\theta:h(\theta)\leq h(\Theta)+\epsilon_n\right\}$, where $\epsilon_n\rightarrow 0$ and
$n\epsilon_n\rightarrow\infty$, as $n\rightarrow\infty$, and proved the following result.

\begin{theorem}
\label{theorem:shalizi3}
Under (S1)--(S7), if for each $\delta>0$, 
\begin{equation}
\tau\left(\mathcal G_n\cap N^c_{\epsilon_n},\delta\right)\leq n
\label{eq:rate0}
\end{equation}
eventually almost surely, then (\ref{eq:conv_rate}) holds almost surely.
\end{theorem}


To investigate the rate of convergence in our cases, 
we need to show that for any $\varepsilon>0$ and all $n$ sufficiently large,
\begin{equation}
\frac{1}{n}\log\int_{\mathcal G_n\cap N^c_{\epsilon_n}}R_n(\theta)d\pi(\theta)\leq -h\left(\mathcal G_n\cap N^c_{\epsilon_n}\right)+\varepsilon.
\label{eq:rate1}
\end{equation}
For $\epsilon_n\downarrow 0$ such that $n\epsilon_n\rightarrow 0$ as $n\rightarrow\infty$, it holds that $N^c_{\epsilon_n}\uparrow \Theta$. Since 
$\mathcal G_n\uparrow \Theta$ as well, $h\left(\mathcal G_n\cap N^c_{\epsilon_n}\right)\downarrow h(\Theta)$, since $h(\theta)$ is continuous in $\theta$.
Combining these arguments with(\ref{eq:rate1}) makes it clear that if we can show
\begin{equation}
\frac{1}{n}\log\int_{\mathcal G_n\cap N^c_{\epsilon_n}}R_n(\theta)d\pi(\theta)\leq -h\left(\Theta\right)+\varepsilon,
\label{eq:rate2}
\end{equation}
for any $\varepsilon>0$ and all $n$ sufficiently large, where $\epsilon_n\downarrow 0$ such that $n\epsilon_n\rightarrow 0$ as $n\rightarrow\infty$,
then that $\epsilon_n$ is the rate of convergence. Now, the same steps as (\ref{eq:s6_2}) lead to 
\begin{align}
\frac{1}{n}\log\int_{\mathcal G_n\cap N^c_{\epsilon_n}}R_n(\btheta)d\pi(\btheta)&\leq\frac{1}{n}\log\left[\underset{\btheta\in\mathcal G_n}{\sup}~R_n(\btheta)\pi(\mathcal G_n)\right]\notag\\
\leq \frac{1}{n}\underset{\btheta\in\mathcal G_n}{\sup}~\log R_n(\btheta). 
\label{eq:rate3}
\end{align}
Since $h(\theta)$ is continuous and coercive for both Gaussian and double exponential errors, in the light of (\ref{eq:rate3}), (\ref{eq:s6_3}), (\ref{eq:s6_4}), (\ref{eq:s6_5}) and (\ref{eq:s6_6}) 
we only need to verify (\ref{eq:s6_7}) to establish (\ref{eq:rate2}). As we have already verified (\ref{eq:s6_7}) for both Gaussian and double exponential errors, (\ref{eq:rate2}) stands verified. 

In other words, (\ref{eq:rate0}), and hence
(\ref{eq:conv_rate}) hold for both the Gaussian process models with Gaussian and double exponential errors, so that their convergence rate is given by $\epsilon_n$.
In other words, the posterior rate of convergence with respect to KL-divergence is just slower than $n^{-1}$ (just slower that $n^{-\frac{1}{2}}$ with respect to
Hellinger distance), for both kinds of errors that we consider.
Our result can be formally stated as the following theorem.
\begin{theorem}
\label{theorem:conv_rate}
For Gaussian process regression with either normal or double exponential errors, under (A1)--(A5),  
(\ref{eq:conv_rate}) holds almost surely, for $\epsilon_n\downarrow 0$ such that $n\epsilon_n\rightarrow\infty$.
\end{theorem}

\subsection{Consequences of model misspecification}
\label{subsec:miss}

Suppose that the true function $\eta_0$ consists of countable number of discontinuities but has continuous first order partial derivatives at all other points. 
Then $\eta_0\notin\mathcal C'(\mathcal X)$, that is, $\eta_0$ is not
in the parameter space. However, there exists some $\tilde\eta\in\mathcal C'(\mathcal X)$ such that $\tilde\eta(\bx)=\eta_0(\bx)$ for all $\bx\in\mathcal X$ where
$\eta_0$ is continuous. Then if the probability measure $Q$ of (A3) is dominated by the Lebesgue measure, it follows from (\ref{eq:h}) and (\ref{eq:h_de}), that $h(\Theta)=0$
for both the Gaussian and double exponential error models. In this case, the posterior of $\eta$ concentrates around $\tilde\eta$, which is the same as $\eta_0$ except at the
countable number of discontinuities of $\eta_0$. If $(\eta_0,\sigma_0)$ is such that $0<h(\Theta)<\infty$, then the posterior concentrates around the minimizers of $h(\theta)$,
provided such minimizers exist in $\Theta$.

Now, following Shalizi, let us define the one-step-ahead
predictive distribution of $\theta$ by $F^n_{\theta}\equiv F_{\theta}\left(Y_n|Y_1,\ldots,Y_{n-1}\right)$, with the convention that $n=1$ 
gives the marginal distribution of the first observation.
Similarly, let $P^n\equiv P^n\left(Y_n|Y_1,\ldots,Y_{n-1}\right)$, which is the best prediction one could make had $P$ been known. The posterior predictive distribution is given by
$F^n_{\pi}=\int_{\Theta}F^n_{\theta}d\pi\left(\theta|\bY_n\right)$.
With the above definitions, \ctn{Shalizi09} proved the following results:
\begin{theorem}
\label{theorem:shalizi4}
Under assumptions (S1)--(S7), with probability 1,
\begin{align}
\underset{n\rightarrow\infty}{\lim\sup}~\rho^2_H\left(P^n,F^n_{\pi}\right)\leq h\left(\Theta\right);\label{eq:hell}\\
\underset{n\rightarrow\infty}{\lim\sup}~\rho^2_{TV}\left(P^n,F^n_{\pi}\right)\leq 4h\left(\Theta\right),\label{eq:tv}
\end{align}
where $\rho_H$ and $\rho_{TV}$ are Hellinger and total variation metrics, respectively.
\end{theorem}
Since, for both our Gaussian process models with normal and double exponential errors, $h(\Theta)=0$ if $\eta_0$ consists of countable number of discontinuities, 
it follows from (\ref{eq:hell}) and (\ref{eq:tv}) that
in spite of such misspecification, the posterior predictive distribution does a good job in learning the best possible predictive distribution
in terms of the popular Hellinger and the total variation distance.
We state our result formally as the following theorem.
\begin{theorem}
\label{theorem:miss}
In the Gaussian process regression problem with either normal or double exponential errors, 
assume that the true function $\eta_0$ consists of countable number of discontinuities but has continuous first order partial derivatives at all other points.
Also assume that the probability measure $Q$ of (A3) is dominated by the Lebesgue measure. Then under (A1) -- (A5), 
\begin{align}
\underset{n\rightarrow\infty}{\lim\sup}~\rho_H\left(P^n,F^n_{\pi}\right)=0;\notag\\
\underset{n\rightarrow\infty}{\lim\sup}~\rho_{TV}\left(P^n,F^n_{\pi}\right)=0,\notag
\end{align}
almost surely.
\end{theorem}

\section{The general nonparametric regression setup}
\label{sec:nongp_setup}
Following \ctn{Choi09} we consider the following model:
\begin{align}
y_i&=\eta(\bx_i)+\epsilon_i;~i=1,\ldots,n;\label{eq:nongp_model1}\\
\epsilon_i &\stackrel{iid}{\sim}\frac{1}{\sigma}\phi\left(\frac{\epsilon_i}{\sigma}\right);~\sigma>0;\label{eq:nongp_model2}\\
\eta(\cdot) &\sim \pi_{\eta}(\cdot);\label{eq:nongp_prior_eta}\\
\sigma &\sim\pi_{\sigma}(\cdot).
\label{eq:nongp_prior_sigma}
\end{align}
In (\ref{eq:nongp_model2}), we model the random errors $\epsilon_i$; $i=1,\ldots,n$ as $iid$ samples from some density $\frac{1}{\sigma}\phi\left(\frac{\cdot}{\sigma}\right)$.
In (\ref{eq:nongp_prior_eta}), $\pi_{\eta}$ stands for any reasonable stochastic process prior, which may may or may not be Gaussian, and in 
(\ref{eq:nongp_prior_sigma}), $\pi_{\sigma}$ is some appropriate prior on $\sigma$.

\subsection{Additional assumptions and their discussions}
\label{subsec:nongp_assumptions}
Regarding the model and the prior, we make the following assumptions in addition to (A1) -- (A5) presented in Section \ref{subsec:assumptions}:
\begin{itemize}

\item[(A6)] The prior on $\eta$ is chosen such that for $\beta>2h\left(\Theta\right)$, 
\begin{align}
	\pi\left(\|\eta\|\leq\exp\left(\left(\beta n\right)^{1/4}\right)\right)&\geq 1-c_{\eta}\exp\left(-\beta n\right);\notag\\
	\pi\left(\|\eta^\prime_j\|\leq\exp\left(\left(\beta n\right)^{1/4}\right)\right)&\geq 1-c_{\eta^\prime_j}\exp\left(-\beta n\right),~\mbox{for}~j=1,\ldots,d,
\end{align}
where $c_{\eta}$ and $c_{\eta^\prime_j}$; $j=1,\ldots,d$, are positive constants.

\item[(A7)] $\phi(\cdot)$ is symmetric about zero; that is, for any $x\in\mathbb R$, $\phi(x)=\phi(|x|)$.  
Further, $\log\phi$ is $L$-Lipschitz, that is, there exists a $L>0$ such that $|\log\phi(x_1)-\log\phi(x_2)|\leq L|x_1-x_2|$, for any $x_1,x_2\in\mathbb R$.

\item[(A8)] For $\bx\in\mathcal X$, let $g_{\eta,\sigma}(\bx)=E_{\theta_0}\left[\log\phi\left(\frac{y-\eta(\bx)}{\sigma}\right)\right]
=\int_{-\infty}^{\infty}\log\phi\left(\frac{\sigma_0z+\eta_0(\bx)-\eta(\bx)}{\sigma}\right)\phi(z)dz$. Then
given $(\eta,\sigma)$, 
$U_i=\log\phi\left(\frac{y_i-\eta(\bx_i)}{\sigma}\right)-g_{\eta,\sigma}(\bx_i)$ are independent sub-exponential random variables satisfying
for any $i=1,\ldots,n$, 
\begin{equation}
E_{\theta_0}\left[\exp\left(\lambda U_i\right)\right]\leq \exp\left(\frac{\lambda^2s^2_{\eta,\sigma}}{2}\right),~\mbox{for}~|\lambda|\leq s^{-1}_{\eta,\sigma}, 
%
\label{eq:a7_1}
\end{equation}
where, for $c_1>0$, $c_2>0$, 
\begin{equation}
s_{\eta,\sigma}= \frac{c_1\|\eta-\eta_0\|+c_2}{\sigma}.
\label{eq:s_eta_sigma}
\end{equation}



\item[(A9)] For $\sigma>0$, $\int_{-\infty}^{\infty}\left|\log\phi\left(\frac{\sigma_0}{\sigma}z\right)\right|\phi(z)dz\leq\frac{c_3}{\sigma}$, where $c_3>0$.
Also, $\int_{-\infty}^{\infty}|z|\phi(z)dz<\infty$.
\item[(A10)]
\begin{enumerate}
\item[(i)]$E_{\bX}\left[g_{\eta,\sigma}(\bX)\right]$ is jointly continuous in $(\eta,\sigma)$;
\item[(ii)]$E_{\bX}\left[g_{\eta,\sigma}(\bX)\right]\rightarrow\infty$ as $\|\theta\|=\|\eta\|+\sigma\rightarrow\infty$.
\end{enumerate}

\end{itemize}

\subsubsection{Discussion of the new assumptions}
\label{subsubsec:nongp_ass_diss}

Condition (A6) ensures that the prior probabilities of the complements of the sieves are exponentially small. Such a requirement is common to most Bayesian asymptotic theories.
In particular, the first two inequalities are satisfied by Gaussian process priors even if $\exp\left(\left(\beta n\right)^{1/4}\right)$ is replaced by $\sqrt{\beta n}$.


Assumption (A7) is the same as that of \ctn{Choi09}, and holds in the case of double exponential errors, for instance.

Conditions (A8), (A9) and (A10) are reasonably mild conditions, and as shown in the supplement, are satisfied by double exponential errors. 

As before, let $\Theta=\mathcal C'(\mathcal X)\times\mathbb R^+$ denote the infinite-dimensional parameter space for our model.

\subsection{Posterior convergence}
\label{subsec:nongp_postconv_normal_de_errors}

As before, we provide a summary of our results leading to posterior convergence in our general setup. The details
are provided in the supplement. 
\begin{lemma}
\label{lemma:nongp_lemma1}
Under our model assumptions and conditions (A1) and (A3), the KL-divergence rate $h(\theta)$ exists for $\theta\in\Theta$, and is given by
\begin{equation}
h(\theta)=\log\left(\frac{\sigma}{\sigma_0}\right)+c-E_\bX\left[g_{\eta,\sigma}(\bX)\right],
\label{eq:nongp_h}
\end{equation}
where $c=E_{\theta_0}\left[\log\phi\left(\frac{y_i-\eta_0(\bx_i)}{\sigma_0}\right)\right]=\int_{-\infty}^{\infty}\left[\log\phi(z)\right]\phi(z)dz$.
\end{lemma}

\begin{theorem}
\label{theorem:nongp_theorem1}
Under our model assumptions and conditions (A1) and (A3), the asymptotic equipartition property holds, and is given by
\begin{equation*}
\underset{n\rightarrow\infty}{\lim}~\frac{1}{n}\log R_n(\theta)=-h(\theta),~\mbox{almost surely}.
\end{equation*}
The convergence is uniform on any compact subset of $\Theta$.
\end{theorem}

Lemma \ref{lemma:nongp_lemma1} and Theorem \ref{theorem:nongp_theorem1} ensure that conditions (S1) -- (S3) of Shalizi hold, and (S4) holds since $h(\theta)$ is almost surely
finite. We construct the sieves $G_n$ as in (\ref{eq:G}). 
Hence, as before, $G_n\rightarrow\Theta$ as $n\rightarrow\infty$ and the assumptions on $\eta$, $\eta'$ given by (A6),
together with (A4) ensure that $\pi(\mathcal G^c_n)\leq\alpha\exp(-\beta n)$, for some $\alpha>0$.
This result, continuity of $h(\theta)$, compactness of $\mathcal G_n$ and the uniform convergence result of Theorem \ref{theorem:theorem1}, together ensure (S5).

As regards (S6), let us note that from the definition of $g_{\eta,\sigma}(\bx)$ and Lipschitz continuity of $\log\phi$, it follows that 
$E_{\bX}\left[g_{\eta,\sigma}(\bX)\right]$ is Lipschitz continuous in $\eta$.
However, we still need to assume that $E_{\bX}\left[g_{\eta,\sigma}(\bX)\right]$ is jointly continuous in $\theta=(\eta,\sigma)$.
Due to (A10) it follows that $h(\theta)$ is continuous in $\theta$ and 
$h(\theta)\rightarrow\infty$ as $\|\theta\|\rightarrow\infty$. 
In other words,
$h(\theta)$ is a continuous coercive function. Hence, $\mathcal S$ is a compact set. 
With these observations, we then have the following result analogous to the Gaussian process case, the proof which is provided in the supplement.
\begin{theorem}
\label{theorem:nongp_theorem3}
In our setup, under (A1)--(A10), it holds that 
\begin{equation*}
\sum_{n=1}^{\infty}\int_{\mathcal S^c}P\left(\left|\frac{1}{n}\log R_n(\theta)+h(\theta)\right|> \kappa-h\left(\Theta\right)\right)d\pi(\theta)<\infty. 
\end{equation*}
\end{theorem}
Since $h(\theta)$ is continuous, (S7) holds trivially. 
Thus, all the assumptions (S1)--(S7) are satisfied, showing that 
Theorems \ref{theorem:shalizi1} and \ref{theorem:shalizi2} hold.
Formally, our results lead to the following theorem.
\begin{theorem}
\label{theorem:nongp2}
Assume the hierarchical model given by (\ref{eq:nongp_model1}), (\ref{eq:nongp_model2}), (\ref{eq:nongp_prior_eta}) and (\ref{eq:nongp_prior_sigma}). 
Then under the conditions (A1) -- (A10), (\ref{eq:post_conv1}) holds. 
Also, for any measurable set $A$ with $\pi(A)>0$, if $\beta>2h(A)$, where $h$ is given by (\ref{eq:nongp_h}), or if $A\subset\cap_{k=n}^{\infty}\mathcal G_k$ for some $n$, 
where $\mathcal G_k$ is given by (\ref{eq:G}), then
(\ref{eq:post_conv2}) holds.
\end{theorem}

\subsection{Rate of convergence and consequences of model misspecification}
\label{subsec:nongp_rate}

For the general nonparametric model, the same result as Theorem \ref{theorem:conv_rate} holds, under (A1)--(A10). Also, the same issues regarding model misspecification
as detailed in Section \ref{subsec:miss} continues to be relevant in this setup. In other words, Theorem \ref{theorem:miss} holds under (A1) -- (A10).

\section{Conclusion}
\label{sec:conclusion}
The fields of both theoretical and applied Bayesian nonparametric regression are dominated by Gaussian process priors and Gaussian noise. 
From the asymptotics perspective, even in the Gaussian setup, a comprehensive theory unifying posterior convergence for both random
and non-random covariates along with the rate of convergence in the case of general priors for the unknown error variance, while also allowing for misspecification, 
seems to be very rare. Even more rare is the aforementioned investigations in the setting where a general stochastic process prior is on the unknown regression function 
is considered and the noise distribution is non-Gaussian and thick-tailed.

The approach of Shalizi allowed us to consider the asymptotic theory incorporating all the above issues,
for both Gaussian and general stochastic process prior for the regression function. The approach,
apart from enabling us to ensure consistency for both random and non-random covariates, allows us to compute the rate of convergence, 
while allowing misspecifications. Perhaps the most interesting result that we obtained is that even if the unknown regression function is misspecified, the posterior
predictive distribution still captures the true predictive distribution asymptotically, for both Gaussian and general setups.

It seems that the most important condition among the assumptions of Shalizi is the asymptotic equipartition property. This directly establishes the KL property 
of the posterior which
characterizes the posterior convergence, the rate of posterior convergence and misspecification.
Interestingly, such a property that plays the key role, turned out to be relatively easy to establish
in our context under reasonably mild conditions. On the other hand, in all the applications
that we investigated so far, (S6) turned out
to be the most difficult to verify. But the approach we devised to handle this condition and
the others, seem to be generally applicable for investigating posterior asymptotics in general
Bayesian parametric and nonparametric problems.

\newpage

\renewcommand\thefigure{S-\arabic{figure}}
\renewcommand\thetable{S-\arabic{table}}
\renewcommand\thesection{S-\arabic{section}}
\renewcommand\thetheorem{S-\arabic{theorem}}

\setcounter{section}{0}
\setcounter{theorem}{0}

\begin{center}
{\LARGE\bf Supplementary Material}
\end{center}

\section{Preliminaries for ensuring posterior consistency under general set-up}
\label{sec:shalizi}

Following \ctn{Shalizi09} we consider a probability space $(\Omega,\mathcal F, P)$, 
and a sequence of random variables $y_1,y_2,\ldots$,   
taking values in some measurable space $(\Xi,\mathcal Y)$, whose
infinite-dimensional distribution is $P$. The natural filtration of this process is
$\sigma(\by_n)$, the smallest $\sigma$-field with respect to which $\by_n$ is measurable.. 

We denote the distributions of processes adapted to $\sigma(\by_n)$ 
by $F_{\theta}$, where $\theta$ is associated with a measurable
space $(\Theta,\mathcal T)$, and is generally infinite-dimensional. 
For the sake of convenience, we assume, as in \ctn{Shalizi09}, that $P$
and all the $F_{\theta}$ are dominated by a common reference measure, with respective
densities $f_{\theta_0}$ and $f_{\theta}$. The usual assumptions that $P\in\Theta$ or even $P$ lies in the support 
of the prior on $\Theta$, are not required for Shalizi's result, rendering it very general indeed.

\subsection{Assumptions and theorems of Shalizi}
\label{subsec:assumptions_shalizi}

\begin{itemize}
\item[(S1)] Consider the following likelihood ratio:
\begin{equation}
R_n(\theta)=\frac{f_{\theta}(\bY_n)}{f_{\theta_0}(\bY_n)}.
\label{eq:R_n}
\end{equation}
Assume that $R_n(\theta)$ is $\sigma(\bY_n)\times \mathcal T$-measurable for all $n>0$.
\end{itemize}

\begin{itemize}
\item[(S2)] For every $\theta\in\Theta$, the KL-divergence rate
\begin{equation}
h(\theta)=\underset{n\rightarrow\infty}{\lim}~\frac{1}{n}E\left(\log\frac{f_{\theta_0}(\bY_n)}{f_{\theta}(\bY_n)}\right).
\label{eq:S3}
\end{equation}
exists (possibly being infinite) and is $\mathcal T$-measurable.
\end{itemize}

\begin{itemize}
\item[(S3)] For each $\theta\in\Theta$, the generalized or relative asymptotic equipartition property holds, and so,
almost surely,
\begin{equation*}
\underset{n\rightarrow\infty}{\lim}~\frac{1}{n}\log R_n(\theta)=-h(\theta).
\end{equation*}
\end{itemize}

\begin{itemize}
\item[(S4)] 
Let $I=\left\{\theta:h(\theta)=\infty\right\}$. 
The prior $\pi$ satisfies $\pi(I)<1$.
\end{itemize}

Following the notation of \ctn{Shalizi09}, for $A\subseteq\Theta$, let
\begin{align}
h\left(A\right)&=\underset{\theta\in A}{\mbox{ess~inf}}~h(\theta);\label{eq:h2}\\
J(\theta)&=h(\theta)-h(\Theta);\label{eq:J}\\
J(A)&=\underset{\theta\in A}{\mbox{ess~inf}}~J(\theta).\label{eq:J2}
\end{align}
\begin{itemize}
\item[(S5)] There exists a sequence of sets $\mathcal G_n\rightarrow\Theta$ as $n\rightarrow\infty$ 
such that: 
\begin{enumerate}
\item[(1)]
\begin{equation}
\pi\left(\mathcal G_n\right)\geq 1-\alpha\exp\left(-\beta n\right),~\mbox{for some}~\alpha>0,~\beta>2h(\Theta);
\label{eq:S5_1}
\end{equation}
\item[(2)]The convergence in (S3) is uniform in $\theta$ over $\mathcal G_n\setminus I$.
\item[(3)] $h\left(\mathcal G_n\right)\rightarrow h\left(\Theta\right)$, as $n\rightarrow\infty$.
\end{enumerate}
\end{itemize}
For each measurable $A\subseteq\Theta$, for every $\delta>0$, there exists a random natural number $\tau(A,\delta)$
such that
\begin{equation}
n^{-1}\log\int_{A}R_n(\theta)\pi(\theta)d\theta
\leq \delta+\underset{n\rightarrow\infty}{\lim\sup}~n^{-1}
\log\int_{A}R_n(\theta)\pi(\theta)d\theta,
\label{eq:limsup_2}
\end{equation}
for all $n>\tau(A,\delta)$, provided 
$\underset{n\rightarrow\infty}{\lim\sup}~n^{-1}\log\pi\left(\mathbb I_A R_n\right)<\infty$.
Regarding this, the following assumption has been made by Shalizi:
\begin{itemize}
\item[(S6)] The sets $\mathcal G_n$ of (S5) can be chosen such that for every $\delta>0$, the inequality
$n>\tau(\mathcal G_n,\delta)$ holds almost surely for all sufficiently large $n$.
\end{itemize}
\begin{itemize}
\item[(S7)] The sets $\mathcal G_n$ of (S5) and (S6) can be chosen such that for any set $A$ with $\pi(A)>0$, 
\begin{equation}
h\left(\mathcal G_n\cap A\right)\rightarrow h\left(A\right),
\label{eq:S7}
\end{equation}
as $n\rightarrow\infty$.
\end{itemize}
Under the above assumptions, \ctn{Shalizi09} proved the following results.

\begin{theorem}[\ctn{Shalizi09}]
\label{theorem:shalizi1}
Consider assumptions (S1)--(S7) and any set $A\in\mathcal T$ with $\pi(A)>0$ and $h(A)>h(\Theta)$. Then,
\begin{equation*}
\underset{n\rightarrow\infty}{\lim}~\pi(A|\bY_n)=0,~\mbox{almost surely.}
\end{equation*}
\end{theorem}

The rate of convergence of the log-posterior is given by the following result. 
\begin{theorem}[\ctn{Shalizi09}]
\label{theorem:shalizi2}
Consider assumptions (S1)--(S7) and any set $A\in\mathcal T$ with $\pi(A)>0$. If $\beta>2h(A)$, where
$\beta$ is given in (\ref{eq:S5_1}) under assumption (S5), or if $A\subset\cap_{k=n}^{\infty}\mathcal G_k$ for some $n$, then
\begin{equation*}
	\underset{n\rightarrow\infty}{\lim}~\frac{1}{n}\log\pi(A|\bY_n)=-J(A),~\mbox{almost surely.}
\end{equation*}
\end{theorem}

\section{Verification of the assumptions of Shalizi for the Gaussian process model with normal errors}
\label{sec:verification}

\subsection{Verification of (S1)}
\label{subsec:S1}
note that
\begin{align}
f_{\theta}(\bY_n)&=\frac{1}{\left(\sigma\sqrt{2\pi}\right)^n}\exp\left\{-\frac{1}{2\sigma^2}\sum_{i=1}^n(Y_i-\eta(\bx_i))^2\right\};\label{eq:like1}\\
f_{\theta_0}(\bY_n)&=\frac{1}{\left(\sigma_0\sqrt{2\pi}\right)^n}\exp\left\{-\frac{1}{2\sigma^2_0}\sum_{i=1}^n(Y_i-\eta_0(\bx_i))^2\right\}.\label{eq:true_like1}
\end{align}
The equations (\ref{eq:like1}) and (\ref{eq:true_like1}) yield, in our case,
\begin{equation}
\frac{1}{n}\log R_n(\theta)=\log\left(\frac{\sigma_0}{\sigma}\right)+\frac{1}{2\sigma^2_0}\times\frac{1}{n}\sum_{i=1}^n\left(y_i-\eta_0(\bx_i)\right)^2
-\frac{1}{2\sigma^2}\times\frac{1}{n}\sum_{i=1}^n\left(y_i-\eta(\bx_i)\right)^2.
\label{eq:R1}
\end{equation}
We show that the right hand side of (\ref{eq:R1}), which we denote as $f(\by_n,\theta)$, is continuous in $(\by_n,\theta)$, which is sufficient to confirm measurability of $R_n(\theta)$.
Let $\|(\by_n,\theta)\|=\|\by_n\|+\|\theta\|$, where $\|\by_n\|$ is the Euclidean norm and $\|\theta\|=\|\eta\|+|\sigma|$, with
$\|\eta\|=\underset{\bx\in\mathcal X}{\sup}~|\eta(\bx)|$. Since $\mathcal X$ is compact and $\eta$ is almost surely continuous, it follows that $\|\eta\|<\infty$ almost surely.

Consider $\by_{n}=(y_{1},y_{2},\ldots,y_{n})^T$ 
and $\boeta_{0n}=(\eta_0(\bx_1),\ldots,\eta_0(\bx_n))^T$. 
Then
\begin{equation}
\sum_{i=1}^n(y_{i}-\eta_0(\bx_i))^2=\by^T_n\by_n-2\by^T_n\boeta_{0n}+\boeta^T_{0n}\boeta_{0n}
\label{eq:cont2}
\end{equation}
is clearly continuous in $\by_n$.
%
Now note that 
\begin{equation}
\frac{1}{n}\sum_{i=1}^n(y_i-\eta(\bx_i))^2=\frac{1}{n}\sum_{i=1}^n(y_i-\eta_0(\bx_i))^2+\frac{1}{n}\sum_{i=1}^n(\eta(\bx_i)-\eta_0(\bx_i))^2
-\frac{2}{n}\sum_{i=1}^n(y_i-\eta_0(\bx_i))(\eta(\bx_i)-\eta_0(\bx_i)),
\label{eq:breakup2}
\end{equation}
where we have already proved continuity of the first term on the right hand side of (\ref{eq:breakup2}). 
To see continuity of $\frac{1}{n}\sum_{i=1}^n(\eta(\bx_i)-\eta_0(\bx_i))^2$ with respect to $\eta$, first consider any sequence $\left\{\eta_j:j=1,2,\ldots\right\}$ satisfying
$\|\eta_j-\tilde\eta\|\rightarrow 0$, as $j\rightarrow\infty$. Then
\begin{align}
&\Bigg|\frac{1}{n}\sum_{i=1}^n(\eta_j(\bx_i)-\eta_0(\bx_i))^2-\frac{1}{n}\sum_{i=1}^n(\tilde\eta(\bx_i)-\eta_0(\bx_i))^2\Bigg |\label{eq:cont3}\\
&\ \ \leq\frac{1}{n}\sum_{i=1}^n\left|\eta_j(\bx_i)-\tilde\eta(\bx_i)\right|\times\left|(\eta_j(\bx_i)-\eta_0(\bx_i))+(\tilde\eta(\bx_i)-\eta_0(\bx_i))\right|\notag\\
&\ \ \leq\|\eta_j-\tilde\eta\|\times\left[\|\eta_j-\eta_0\|+\|\tilde\eta-\eta_0\|\right]\notag\\
&\ \ \leq\|\eta_j-\tilde\eta\|\times\left[\|\eta_j-\tilde\eta\|+2\|\tilde\eta-\eta_0\|\right]\notag\\
&\ \ \rightarrow 0,~\mbox{as}~j\rightarrow\infty.
\label{eq:cont4}
\end{align}
This proves continuity of the second term of (\ref{eq:breakup2}). 

For the third term of (\ref{eq:breakup2}) we now prove that for any $\tilde \by\in\mathbb R^n$, and for any sequence 
$\left\{\by_j:j=1,2,\ldots\right\}$ (we denote the $i$-th component of $\by_j$ as $y_{ij}$) such that
$\|\by_j-\tilde\by\|\rightarrow 0$, as $j\rightarrow\infty$, and for any function $\tilde\eta$ associated with any sequence $\left\{\eta_j:j=1,2,\ldots\right\}$ satisfying
$\|\eta_j-\tilde\eta\|\rightarrow 0$, as $j\rightarrow\infty$, $\sum_{i=1}^n(y_{ij}-\eta_0(\bx_i))(\eta_j(\bx_i)-\eta_0(\bx_i))\rightarrow
\sum_{i=1}^n(\tilde y_i-\eta_0(\bx_i))(\tilde\eta(\bx_i)-\eta_0(\bx_i))$, as $j\rightarrow\infty$.
Indeed, observe that
\begin{align}
&\Bigg|\sum_{i=1}^n(y_{ij}-\eta_0(\bx_i))(\eta_j(\bx_i)-\eta_0(\bx_i))-\sum_{i=1}^n(\tilde y_i-\eta_0(\bx_i))(\tilde\eta(\bx_i)-\eta_0(\bx_i))\Bigg|\notag\\
&\qquad =\Bigg|\sum_{i=1}^n(y_{ij}-\tilde y_i)(\eta_j(\bx_i)-\tilde\eta(\bx_i))
+\sum_{i=1}^n(\tilde y_i-\eta_0(\bx_i))(\eta_j(\bx_i)-\tilde\eta(\bx_i))\notag\\
&\qquad\qquad +\sum_{i=1}^n(y_{ij}-\tilde y_i)(\tilde\eta(\bx_i)-\eta_0(\bx_i))\Bigg|\notag\\
&\qquad\leq n\|\by_j-\tilde\by\|\|\eta_j-\tilde\eta\|+
n\|\tilde\by-\boeta_{0n}\|\|\eta_j-\tilde\eta\|+n\|\by_j-\tilde\by\|\|\tilde\eta-\eta_0\|\notag\\
&\qquad\rightarrow 0,~\mbox{as}~\|\by_j-\tilde\by\|\rightarrow 0~\mbox{and}~\|\eta_j-\tilde\eta\|\rightarrow 0,~\mbox{as}~j\rightarrow\infty.\notag
\end{align}
Hence, $\sum_{i=1}^n(y_i-\eta_0(\bx_i))(\eta(\bx_i)-\eta_0(\bx_i))$ is continuous in $\by_n$ and $\eta$.
Continuity is clearly preserved if the above expression is divided by $\sigma$.


Also, the first term of $f(\by_n,\theta)$, given by $\log\left(\frac{\sigma_0}{\sigma}\right)$, is clearly continuous
in $\sigma$. 
Thus, continuity of $f(\by_n,\theta)$ with respect to $(\by_n,\theta)$ is guaranteed, so that (S1) holds.
Also observe that when the covariates are regarded as random, due to measurability of $\eta_0(\bX)$ assumed in (A4) and continuity of $\eta(\bx)$ in $\bx$.

\subsection{Verification of (S2) and proof of Lemma \ref{lemma:lemma1} for Gaussian errors}
\label{subsec:S2}
It follows from (\ref{eq:like1}) and (\ref{eq:true_like1}), that
\begin{equation}
\log\frac{f_{\theta_0}(\by_n)}{f_{\theta}(\by_n)}=n\log\left(\frac{\sigma}{\sigma_0}\right)-\frac{1}{2\sigma^2_0}\sum_{i=1}^n(y_i-\eta_0(\bx_i))^2+
\frac{1}{2\sigma^2}\sum_{i=1}^n(y_i-\eta(\bx_i))^2,
\label{eq:logratio1}
\end{equation}
so that
\begin{equation}
\frac{1}{n}E_{\theta_0}\left(\log\frac{f_{\theta_0}(\by_n)}{f_{\theta}(\by_n)}\right)=\log\left(\frac{\sigma}{\sigma_0}\right)-\frac{1}{2}
+\frac{\sigma^2_0}{2\sigma^2}+\frac{1}{2\sigma^2}\times\frac{1}{n}\sum_{i=1}^n\left(\eta(\bx_i)-\eta_0(\bx_i)\right)^2.
\label{eq:logratio2}
\end{equation}
By (A3), as $n\rightarrow\infty$,
\begin{equation}
\frac{1}{n}\sum_{i=1}^n\left(\eta(\bx_i)-\eta_0(\bx_i)\right)^2\rightarrow E_\bX\left[\eta(\bX)-\eta_0(\bX)\right]^2
=\int_{\mathcal X}\left[\eta(\bX)-\eta_0(\bX)\right]^2dQ.
\label{eq:slln1}
\end{equation}
Hence,
\begin{equation}
\frac{1}{n}E_{\theta_0}\left(\log\frac{f_{\theta_0}(\by_n)}{f_{\theta}(\by_n)}\right)\rightarrow \log\left(\frac{\sigma}{\sigma_0}\right)-\frac{1}{2}
+\frac{\sigma^2_0}{2\sigma^2}+\frac{1}{2\sigma^2}E_\bX\left[\eta(\bX)-\eta_0(\bX)\right]^2,~\mbox{as}~n\rightarrow\infty.
\label{eq:logratio3}
\end{equation}
We let
\begin{equation*}
h(\theta)=\log\left(\frac{\sigma}{\sigma_0}\right)-\frac{1}{2}+\frac{\sigma^2_0}{2\sigma^2}+\frac{1}{2\sigma^2}E_\bX\left[\eta(\bX)-\eta_0(\bX)\right]^2.
\end{equation*}

\subsection{Verification of (S3) and proof of Theorem \ref{theorem:theorem1} for Gaussian errors}
\label{subsec:S3}

By SLLN, as $n\rightarrow\infty$,
\begin{equation}
\frac{1}{n}\sum_{i=1}^n\left(y_i-\eta_0(\bx_i)\right)^2\stackrel{a.s.}{\longrightarrow}\sigma^2_0,
\label{eq:slln2}
\end{equation}
where $``\stackrel{a.s.}{\longrightarrow}"$ denotes convergence almost surely.
Also,
\begin{align}
\frac{1}{n}\sum_{i=1}^n\left(y_i-\eta(\bx_i)\right)^2 &=\frac{1}{n}\sum_{i=1}^n\left(y_i-\eta_0(\bx_i)\right)^2+\frac{1}{n}\sum_{i=1}^n\left(\eta(\bx_i)-\eta_0(\bx_i)\right)^2\notag\\
&\qquad\qquad+\frac{2}{n}\sum_{i=1}^n\left(y_i-\eta_0(\bx_i)\right)\left(\eta_0(\bx_i)-\eta(\bx_i)\right).
\label{eq:breakup1}
\end{align}
By (\ref{eq:slln2}) the first term on the right hand side of (\ref{eq:breakup1}) converges almost surely to $\sigma^2_0$. The second term converges to 
$E_{\bX}\left[\eta(\bX)-\eta_0(\bX)\right]^2$ and the
third term converges almost surely to zero by Kolmogorov's SLLN for independent random variables, noting that $y_i-\eta_0(\bx_i)=\e_i$ are independent zero mean random variables
and $\sum_{i=1}^{\infty}i^{-2}Var\left((y_i-\eta_0(\bx_i)(\eta_0(\bx_i)-\eta(\bx_i))\right)=\sigma^2_0\sum_{i=1}^{\infty}i^{-2}\left(\eta_0(\bx_i)-\eta(\bx_i)\right)^2\leq
\sigma^2_0\|\eta-\eta_0\|^2\sum_{i=1}^{\infty}i^{-2}<\infty$. 
Hence, letting $n\rightarrow\infty$ in (\ref{eq:R1}), it follows that
\begin{equation}
\frac{1}{n}\log R_n(\theta)\stackrel{a.s.}{\longrightarrow}\log\left(\frac{\sigma_0}{\sigma}\right)+\frac{1}{2}-\frac{\sigma^2_0}{2\sigma^2}
-\frac{1}{2\sigma^2}E_{\bX}\left[\eta(\bX)-\eta_0(\bX)\right]^2=-h(\theta).
\label{eq:R2}
\end{equation}
The above results of course remain the same if the covariates are assumed to be random.

\subsection{Verification of (S4)}
\label{subsec:S4}

Note that $h(\theta)\leq\log\left(\frac{\sigma}{\sigma_0}\right)-\frac{1}{2}+\frac{\sigma^2_0}{2\sigma^2}+\frac{\|\eta-\eta_0\|^2}{2\sigma^2}$, where
$0<\sigma<\infty$ and $0<\|\eta-\eta_0\|<\infty$ with prior probability one. Hence, $h(\theta)<\infty$ with probability one, showing that (S4) holds.

\subsection{Verification of (S5)}
\label{subsec:S5}

\subsubsection{Verification of (S5) (1)}
\label{subsubsec:S5_1}


Recall that
\begin{align}
	\mathcal G_n&=\left\{\left(\eta,\sigma\right):\|\eta\|\leq\exp(\left(\beta n\right)^{1/4}),\exp(-\left(\beta n\right)^{1/4})\leq\sigma\leq\exp(\left(\beta n\right)^{1/4}),
	\|\eta'_j\|\leq\exp(\left(\beta n\right)^{1/4});j=1,\ldots,d\right\}.\notag
\end{align}
Then $\mathcal G_n\rightarrow\Theta$, as $n\rightarrow\infty$.
Now note that
\begin{align}
	&\pi(\mathcal G_n)=\pi\left(\|\eta\|\leq\exp(\left(\beta n\right)^{1/4}),\exp(-\left(\beta n\right)^{1/4})\leq\sigma\leq\exp(\left(\beta n\right)^{1/4})\right)\notag\\
	&\ \ -\pi\left(\left\{\|\eta'_j\|\leq\exp(\left(\beta n\right)^{1/4});j=1,\ldots,d\right\}^c\right)\notag\\
	&=\pi\left(\|\eta\|\leq\exp(\left(\beta n\right)^{1/4}),\exp(-\left(\beta n\right)^{1/4})\leq\sigma\leq\exp(\left(\beta n\right)^{1/4})\right)\notag\\
	&\ \ -\pi\left(\bigcup_{j=1}^d\left\{\|\eta'_j\|>\exp(\left(\beta n\right)^{1/4})\right\}\right)\notag\\
	& \geq 1-\pi\left(\|\eta\|>\exp(\left(\beta n\right)^{1/4})\right)-\pi\left(\left\{\exp(-\left(\beta n\right)^{1/4})\leq\sigma\leq\exp(\left(\beta n\right)^{1/4})\right\}^c\right)\notag\\
	&\ \ \  \ -\sum_{j=1}^d\pi\left(\|\eta'_j\|>\exp(\left(\beta n\right)^{1/4})\right)\notag\\
&\geq 1-(c_{\eta}+c_{\sigma}+\sum_{j=1}^dc_{\eta^\prime_j})\exp(-\beta n),
\label{eq:s5_1}
\end{align}
by the Borell-TIS inequality and (A5). 
In other words, (S5) (1) holds.

\subsubsection{Verification of (S5) (2)}
\label{subsubsec:S5_2}

We now show that (S5) (2), namely, convergence in (S3) is uniform in $\theta$ over $\mathcal G_n\setminus I$ holds. First note that $I=\emptyset$ in our case, so that
$\mathcal G_n\setminus I=\mathcal G_n$. 

To proceed further, we show that $\mathcal G_n$ is compact.
Note that $\mathcal G_n=\mathcal G_{n,\eta}\times \mathcal G_{n,\sigma}$, where
\begin{align}
	\mathcal G_{n,\eta}&=\left\{\eta:\|\eta\|\leq\exp(\left(\beta n\right)^{1/4}),~\|\eta'_j\|\leq\exp(\left(\beta n\right)^{1/4});j=1,\ldots,d\right\}\notag
\end{align}
and
\begin{equation*}
	\mathcal G_{n,\sigma}=\left\{\sigma:\exp(-\left(\beta n\right)^{1/4})\leq\sigma\leq\exp(\left(\beta n\right)^{1/4})\right\}.
\end{equation*}
Since $\mathcal G_{n,\sigma}$ is compact and products of compact sets is compact, it is enough to prove compactness of $\mathcal G_{n,\eta}$.
We use the Arzela-Ascoli lemma to prove that $\mathcal G_{n,\eta}$ is compact for each $n\geq 1$. 
In other words, $\mathcal G_{n,\eta}$ is compact if and only if it is closed,
bounded and equicontinuous. By boundedness we mean $|\eta(\bx)|<M$ for each $\bx\in\mathcal X$ and for each $\eta\in\mathcal G_{n,\eta}$. Equicontinuity entails
that for any $\epsilon>0$, there exists $\delta>0$ which depends only on $\epsilon$ such that $|\eta(\bx_1)-\eta(\bx_2)|<\epsilon$ whenever $\|\bx_1-\bx_2\|<\delta$,
for all $\eta\in\mathcal G_{n,\eta}$. Closedness and boundedness are obvious from the definition of $\mathcal G_{n,\eta}$. Equicontinuity follows from the fact that
the elements of $\mathcal G_{n,\eta}$ are Lipschitz continuous thanks to boundedness of the partial derivatives. Thus, $\mathcal G_{n,\eta}$, and hence $\mathcal G_n$ is compact.

Since $\mathcal G_n$ is compact for all $n\geq 1$, uniform convergence as required will be proven if we can show that $\frac{1}{n}\log R_n(\theta)+h(\theta)$ is 
stochastically equicontinuous almost surely in $\theta\in\mathcal G$ for any 
$\mathcal G\in\left\{\mathcal G_n:n=1,2,\ldots\right\}$ and $\frac{1}{n}\log R_n(\theta)+h(\theta)\rightarrow 0$, almost surely, for all $\theta\in\mathcal G$
(see \ctn{Newey91}, \ctn{Billingsley13}) for the general theory of uniform convergence in compact sets under stochastic equicontinuity). 
Since, in the context of (S3) we have already shown almost sure pointwise convergence of $\frac{1}{n}\log R_n(\theta)$ to $-h(\theta)$, it is enough to verify stochastic
equicontinuity of $\frac{1}{n}\log R_n(\theta)+h(\theta)$ in $\mathcal G\in\left\{\mathcal G_n:n=1,2,\ldots\right\}$.

Stochastic equicontinuity usually follows easily if one can prove that the function concerned is almost surely Lipschitz continuous.
Recall
from (\ref{eq:R1}), (\ref{eq:cont2}), (\ref{eq:breakup2}) and (\ref{eq:cont4}) that if the term $\frac{1}{n}\sum_{i=1}^n(y_i-\eta_0(\bx_i))(\eta(\bx_i)-\eta_0(\bx_i))$
can be proved Lipschitz continuous in $\eta\in\mathcal G$, then $\frac{1}{n}\log R_n(\theta)$ is Lipschitz for $\eta\in\mathcal G$. Also, if $E_{\bX}\left[\eta(\bX)-\eta_0(\bX)\right]^2$
is Lipschitz in $\eta$, then it would follow from (\ref{eq:h}) that $h(\theta)$ is Lipschitz for $\eta\in\mathcal G$. Since sum of Lipschitz functions is Lipschitz, this would imply
that $\frac{1}{n}\log R_n(\theta)+h(\theta)$ is Lipschitz in $\eta\in\mathcal G$. Since the first derivative of $\frac{1}{n}\log R_n(\theta)+h(\theta)$
with respect to $\sigma$ is bounded (as $\sigma$ is bounded in $\mathcal G$), it would then follow that $\frac{1}{n}\log R_n(\theta)+h(\theta)$ is Lipschitz for $\theta\in\mathcal G$. 
Hence, to see that $\frac{1}{n}\sum_{i=1}^n(y_i-\eta_0(\bx_i))(\eta(\bx_i)-\eta_0(\bx_i))$ is almost surely Lipschitz in $\eta\in\mathcal G$, note that for any $\eta_1,\eta_2\in\mathcal G$,
\begin{align}
&\left|\frac{1}{n}\sum_{i=1}^n(y_i-\eta_0(\bx_i))(\eta_1(\bx_i)-\eta_0(\bx_i))-\frac{1}{n}\sum_{i=1}^n(y_i-\eta_0(\bx_i))(\eta_2(\bx_i)-\eta_0(\bx_i))\right|\notag\\
&\qquad\qquad\leq\|\eta_1-\eta_2\|\times\frac{1}{n}\sum_{i=1}^n\left|y_i-\eta_0(\bx_i)\right|.\notag
\end{align}
Hence, $\frac{1}{n}\sum_{i=1}^n(y_i-\eta_0(\bx_i))(\eta(\bx_i)-\eta_0(\bx_i)$ is Lipschitz in $\eta$ and since $\frac{1}{n}\sum_{i=1}^n\left|y_i-\eta_0(\bx_i)\right|\rightarrow
E_{\theta_0}\left|y_1-\eta_0(\bx_1)\right| 
<\infty$ as $n\rightarrow\infty$, stochastic equicontinuity follows.

That $E_{\bX}\left[\eta(\bX)-\eta_0(\bX)\right]^2$ is also Lipschitz in $\mathcal G$ can be seen from the fact that for $\eta_1,\eta_2\in\mathcal G$,
\begin{equation*}
\left|E_{\bX}\left[\eta_1(\bX)-\eta_0(\bX)\right]^2-E_{\bX}\left[\eta_2(\bX)-\eta_0(\bX)\right]^2\right|
\leq\|\eta_1-\eta_2\|\times\left[\|\eta_1\|+\|\eta_2\|+2\|\eta_0\|\right],
\end{equation*}
where $\|\eta_0\|<\kappa_0$ by (A4) and for $j=1,2$, 
$\|\eta_j\|\leq \exp(\left(\beta m\right)^{1/4})$,
where $\mathcal G=\mathcal G_m$, for $m\geq 1$.

\subsubsection{Verification of (S5) (3)}
\label{subsubsec:S5_3}

We now verify (S5) (3). For our purpose, let us show that $h(\theta)$ is continuous in $\theta$. Continuity will easily follow if we can show that $E_{\bX}\left[\eta(\bX)-\eta_0(\bX)\right]^2$
is continuous in $\eta$. As before, let $\eta_j$ be a sequence of functions converging to $\tilde\eta$ in the sense $\|\eta_j-\tilde\eta\|\rightarrow 0$ as $j\rightarrow\infty$.
Then, since $\left|E_{\bX}\left[\eta_j(\bX)-\eta_0(\bX)\right]^2-E_{\bX}\left[\tilde\eta(\bX)-\eta_0(\bX)\right]^2\right|
\leq \|\eta_j-\tilde\eta\|\left[\|\eta_j-\tilde\eta\|+2\|\tilde\eta-\eta_0\|\right]
\rightarrow 0$ as $j\rightarrow\infty$, continuity follows. Hence, continuity of $h(\theta)$, compactness of $\mathcal G_n$, along with its non-decreasing nature with respect to $n$ implies that 
$h\left(\mathcal G_n\right)\rightarrow h(\Theta)$, as $n\rightarrow\infty$.

\subsection{Verification of (S6) and proof of Theorem \ref{theorem:theorem3} for Gaussian errors}
\label{subsec:S6}

Observe that 
\begin{align}
\frac{1}{n}\log R_n(\theta)+h(\theta)&=\left[\frac{1}{2\sigma^2_0}\times\frac{1}{n}\sum_{i=1}^n\left(y_i-\eta_0(\bx_i)\right)^2-\frac{1}{2}\right]
+\left[\frac{1}{2\sigma^2}\times\frac{1}{n}\sum_{i=1}^n\left(y_i-\eta_0(\bx_i)\right)^2-\frac{\sigma^2_0}{2\sigma^2}\right]\notag\\
&\qquad+\left[\frac{1}{2\sigma^2}\times\frac{1}{n}\sum_{i=1}^n\left(\eta(\bx_i)-\eta_0(\bx_i)\right)^2-\frac{1}{2\sigma^2}E_{\bX}\left(\eta(\bX)-\eta_0(\bX)\right)^2\right]\notag\\
&\qquad+\left[\frac{1}{\sigma^2}\times\frac{1}{n}\sum_{i=1}^n\left(y_i-\eta_0(\bx_i)\right)\left(\eta(\bx_i)-\eta_0(\bx_i)\right)\right].
\label{eq:s5_5}
\end{align}
Let $\kappa_1=\kappa-h(\Theta)$. Then it follows from (\ref{eq:s5_5}) that for all $\theta\in\mathcal G$, we have
\begin{align}
&P\left(\left|\frac{1}{n}\log R_n(\theta)+h(\theta)\right|>\kappa_1\right)\notag\\
&\leq P\left(\left|\frac{1}{2\sigma^2_0}\times\frac{1}{n}\sum_{i=1}^n\left(y_i-\eta_0(\bx_i)\right)^2-\frac{1}{2}\right|>\frac{\kappa_1}{4}\right)
+P\left(\left|\frac{1}{2\sigma^2}\times\frac{1}{n}\sum_{i=1}^n\left(y_i-\eta_0(\bx_i)\right)^2-\frac{\sigma^2_0}{2\sigma^2}\right|>\frac{\kappa_1}{4}\right)\notag\\
&\qquad+P\left(\left|\frac{1}{2\sigma^2}\times\frac{1}{n}\sum_{i=1}^n\left(\eta(\bx_i)-\eta_0(\bx_i)\right)^2
-\frac{1}{2\sigma^2}E_{\bX}\left(\eta(\bX)-\eta_0(\bX)\right)^2\right|>\frac{\kappa_1}{4}\right)\notag\\
&\qquad+P\left(\left|\frac{1}{\sigma^2}\times\frac{1}{n}\sum_{i=1}^n\left(y_i-\eta_0(\bx_i)\right)\left(\eta(\bx_i)-\eta_0(\bx_i)\right)\right|>\frac{\kappa_1}{4}\right).
\label{eq:s5_6}
\end{align}

Note that $\sum_{i=1}^n\left(\frac{y_i-\eta_0(\bx_i)}{\sigma_0}\right)^2=\bz^T_n\bz_n$, 
where $\bz_n\sim N_n\left(\bzero_n,\bI_n\right)$, the $n$-dimensional normal distribution with mean $\bzero_n=(0,0,\ldots,0)^T$ and covariance matrix $\bI_n$, the identity matrix. 
Using the Hanson-Wright inequality we bound the
first term of the right hand side of (\ref{eq:s5_6}) as follows:
\begin{align}
&P\left(\left|\frac{1}{2\sigma^2_0}\times\frac{1}{n}\sum_{i=1}^n\left(y_i-\eta_0(\bx_i)\right)^2-\frac{1}{2}\right|>\frac{\kappa_1}{4}\right)\notag\\
&\qquad=P\left(\left|\bz^T_n\bz_n-n\right|>\frac{n\kappa_1}{2}\right)\notag\\
&\qquad\leq 2\exp\left(-n\min\left\{\frac{\kappa^2_1}{16c_0},\frac{\kappa_1}{4c_0}\right\}\right),
\label{eq:s5_7}
\end{align}
where $c_0>0$ is a constant.
It follows from (\ref{eq:s5_7}) that
\begin{equation}
\int_{\mathcal S^c}P\left(\left|\frac{1}{2\sigma^2_0}\times\frac{1}{n}\sum_{i=1}^n\left(y_i-\eta_0(\bx_i)\right)^2-\frac{1}{2}\right|>\frac{\kappa_1}{4}\right)d\pi(\theta)
\leq 2\exp\left(-n\min\left\{\frac{\kappa^2_1}{16c_0},\frac{\kappa_1}{4c_0}\right\}\right).
\label{eq:s5_7_int}
\end{equation}

In almost the same way as in (\ref{eq:s5_7}), the second term of the right hand side of (\ref{eq:s5_6}) can be bounded as:
\begin{align}
&P\left(\left|\frac{1}{2\sigma^2}\times\frac{1}{n}\sum_{i=1}^n\left(y_i-\eta_0(\bx_i)\right)^2-\frac{\sigma^2_0}{2\sigma^2}\right|>\frac{\kappa_1}{4}\right)\notag\\
&\qquad=P\left(\left|\bz^T_n\bz_n-n\right|>\frac{n\kappa_1\sigma^2}{2\sigma^2_0}\right)\notag\\
&\qquad\leq 2\exp\left(-n\min\left\{\frac{\kappa^2_1\sigma^4}{16c_0\sigma^4_0},\frac{\kappa_1\sigma^2}{4c_0\sigma^2_0}\right\}\right).
\label{eq:s5_8}
\end{align}
Now
\begin{align}
&\int_{\mathcal S^c}P\left(\left|\frac{1}{2\sigma^2}\times\frac{1}{n}\sum_{i=1}^n\left(y_i-\eta_0(\bx_i)\right)^2-\frac{\sigma^2_0}{2\sigma^2}\right|>\frac{\kappa_1}{4}\right)d\pi(\theta)\notag\\
&\leq \int_{\mathcal G_n}2\exp\left(-n\min\left\{\frac{\kappa^2_1\sigma^4}{16c_0\sigma^4_0},\frac{\kappa_1\sigma^2}{4c_0\sigma^2_0}\right\}\right)\pi(\sigma^2)d\sigma^2\notag\\
&\qquad+\int_{\mathcal G^c_n}2\exp\left(-n\min\left\{\frac{\kappa^2_1\sigma^4}{16c_0\sigma^4_0},\frac{\kappa_1\sigma^2}{4c_0\sigma^2_0}\right\}\right)\pi(\theta)d\theta\notag\\
&\leq \int_{\mathcal G_n}2\exp\left(-n\min\left\{\frac{\kappa^2_1\sigma^4}{16c_0\sigma^4_0},\frac{\kappa_1\sigma^2}{4c_0\sigma^2_0}\right\}\right)\pi(\sigma^2)d\sigma^2
+2\pi(\mathcal G^c_n)\notag\\
	&\leq\int_{\exp(-2\left(\beta n\right)^{1/4})}^{\exp(2\left(\beta n\right)^{1/4})}2\exp\left(-n\frac{\kappa^2_1\sigma^4}{16c_0\sigma^4_0}\right)\pi(\sigma^2)d\sigma^2\notag\\
	&\qquad+\int_{\exp(-2\left(\beta n\right)^{1/4})}^{\exp(2\left(\beta n\right)^{1/4})}2\exp\left(-n\frac{\kappa_1\sigma^2}{4c_0\sigma^2_0}\right)\pi(\sigma^2)d\sigma^2+2\pi(\mathcal G^c_n)\notag\\
	&=\int_{\exp(-2\left(\beta n\right)^{1/4})}^{\exp(2\left(\beta n\right)^{1/4})}2\exp\left(-n\frac{\kappa^2_1u^{-2}}{16c_0\sigma^4_0}\right)\pi(u^{-1})u^{-2}du\notag\\
	&\qquad+\int_{\exp(-2\left(\beta n\right)^{14})}^{\exp(2\left(\beta n\right)^{1/4})}2\exp\left(-n\frac{\kappa_1u^{-1}}{4c_0\sigma^2_0}\right)\pi(u^{-1})u^{-2}du+2\pi(\mathcal G^c_n).
\label{eq:s5_8_int}
\end{align}

Let us first consider the first term of (\ref{eq:s5_8_int}). Note that the prior $\pi\left(u^{-1}\right)u^{-2}$ is such that large values of $u$ receive
small probabilities. Hence, if this prior is replaced by an appropriate function which has a thicker tail than the prior, then the resultant integral provides an upper bound
for the first term of (\ref{eq:s5_8_int}). We consider a function 
$\tilde\pi(u)$ which is of mixture form depending upon $n$, that is, 
we let $\tilde\pi_n(u)=c_3\sum_{r=1}^{M_n}\psi^{\zeta_{rn}}_{rn}\exp(-\psi_{rn} u^2)u^{2(\zeta_{rn}-1)}$, 
where 
$M_n\leq\exp(\left(\beta n\right)^{1/4})$ is the number of mixture components, $c_3>0$,  
for $r=1,\ldots,M_n$, $\frac{1}{2}<\zeta_{rn}\leq c_4n^q$, for $0<q<1/4$ and $n\geq 1$, where $c_4>0$, and $0<\psi_1\leq\psi_{rn}<c_5<\infty$, for all $r$ and $n$. 
In this case, with $C_1=\frac{1}{16c_0\sigma^4_0}$,
\begin{align}
	& \int_{\exp(-2\left(\beta n\right)^{1/4})}^{\exp(2\left(\beta n\right)^{1/4}}\exp\left(-C_1\kappa^2_1nu^{-2}\right)\pi(u^{-1})u^{-2}du\notag\\
&\leq c_3\sum_{r=1}^{M_n}\psi^{\zeta_{rn}}_{rn}
	\int_{\exp(-2\left(\beta n\right)^{1/4})}^{\exp(2\left(\beta n\right)^{1/4})}\exp\left[-\left(C_1\kappa^2_1nu^{-2}+\psi_1u^2\right)\right]\left(u^2\right)^{\zeta_{rn}-1}du.
\label{eq:term1_bound3_1}
\end{align}
Now the $r$-th integrand of (\ref{eq:term1_bound3_1}) is maximized at 
$\tilde u^2_{rn}= \frac{\zeta_{rn}-1+\sqrt{(\zeta_{rn}-1)^2+4C_1\psi_{1}\kappa^2_1 n}}{2\psi_{1}}$, so that for sufficiently large $n$, 
$c_1\kappa_1\sqrt{\frac{n}{\psi_{1}}}\leq\tilde u^2_{rn}\leq \tilde c_1\kappa_1\sqrt{\frac{n}{\psi_{1}}}$, for some positive constants $c_1$ and $\tilde c_1$. 
Now, for sufficiently large $n$, we have $\frac{\tilde u^2_{rn}}{\log\tilde u^2_{rn}}\geq\frac{\zeta_{rn}-1}{\psi_{1}(1-c_2)}$, for $0<c_2<1$.
Hence, for sufficiently large $n$, $C_1\kappa^2_1n\tilde u^{-2}_{rn}+\psi_{1}\tilde u^2_{rn}-(\zeta_{rn}-1)\log(\tilde u^2_{rn})\geq c_2\psi_1\tilde u^2_{rn}
\geq C_2\kappa_1\sqrt{\psi_{1} n}$ for some 
positive constant $C_2$. From these and (\ref{eq:term1_bound3_1}) it follows that
\begin{align}
	&\int_{\exp(-2\left(\beta n\right)^{1/4})}^{\exp(2\left(\beta n\right)^{1/4})}2\exp\left(-n\frac{\kappa^2_1u^{-2}}{16c_0\sigma^4_0}\right)\pi(u^{-1})u^{-2}du\notag\\
&=c_3\sum_{r=1}^{M_n}\psi^{\zeta_{rn}}_{rn}
	\int_{\exp(-2\left(\beta n\right)^{1/4})}^{\exp(2\left(\beta n\right)^{1/4})}\exp\left[-\left(C_1\kappa^2_1nu^{-2}+\psi_1u^2\right)\right]\left(u^2\right)^{\zeta_{rn}-1}du\notag\\
	&\leq c_3M_n\exp\left[-\left(C_2\kappa_1\sqrt{n\psi_1}-2\left(\beta n\right)^{1/4}-\tilde c_5 n^q\right)\right]\notag\\
	&\leq c_3\exp\left[-\left(C_2\kappa_1\sqrt{n\psi_1}-3\left(\beta n\right)^{1/4}-\tilde c_5 n^q\right)\right].
\label{eq:term1_bound3_2}
\end{align}
for some constant $\tilde c_5$. The negative of the exponent of (\ref{eq:term1_bound3_2}) is clearly positive for large $n$. 

For the second term of (\ref{eq:s5_8_int}), we consider $\tilde\pi_n(u)=c_3\sum_{r=1}^{M_n}\psi^{\zeta_{rn}}_{rn}\exp(-\psi_{rn} u)u^{(\zeta_{rn}-1)}$, with
$M_n\leq\exp(\left(\beta n\right)^{1/4})$ being the number of mixture components, $c_3>0$, for $r=1,\ldots,M_n$, $0<\zeta_{rn}\leq c_4n^q$, for $0<q<1/4$ and $n\geq 1$, 
where $c_4>0$, and $0<\psi_1\leq\psi_{rn}<c_5<\infty$, for all $r$ and $n$. Thus, the only difference here with the previous definition of $\tilde\pi_n(u)$
is that here $\zeta_{rn}>0$ instead of $\zeta_{rn}>\frac{1}{2}$, which is due to the fact that here $u^2$ is replaced with $u$. In the same way as in (\ref{eq:term1_bound3_2}), it then
follows that
\begin{equation}
	\int_{\exp(-2\left(\beta n\right)^{1/4})}^{\exp(2\left(\beta n\right)^{1/4})}2\exp\left(-n\frac{\kappa_1u^{-1}}{4c_0\sigma^2_0}\right)\pi(u^{-1})u^{-2}du
	\leq c_3\exp\left[-\left(C_2\sqrt{\kappa_1n\psi_1}-3\left(\beta n\right)^{1/4}-\tilde c_5 n^q\right)\right].
\label{eq:term1_bound3_3}
\end{equation}
Again, the negative of the exponent of (\ref{eq:term1_bound3_3}) is clearly positive for large $n$. 

For the third term, let us first consider the case of random covariates $\bX$. Here observe that by Hoeffding's inequality (\ctn{Hoeffding63}),
\begin{align}
&P\left(\left|\frac{1}{2\sigma^2}\times\frac{1}{n}\sum_{i=1}^n\left(\eta(\bx_i)-\eta_0(\bx_i)\right)^2
-\frac{1}{2\sigma^2}E_{\bX}\left(\eta(\bX)-\eta_0(\bX)\right)^2\right|>\frac{\kappa_1}{4}\right)\notag\\
&\qquad\leq 2\exp\left\{-\frac{nC\kappa^2_1\sigma^4}{\|\eta-\eta_0\|^2}\right\},
\label{eq:zero_prob2}
\end{align}
where $C>0$ is a constant. Note that $\|\eta-\eta_0\|$ is clearly the upper bound of $|\eta(\cdot)-\eta_0(\cdot)|$. 
Such an upper bound is finite since $\mathcal X$ is compact, $\eta(\cdot)$ is continuous on $\mathcal X$,
and $\|\eta_0\|<\infty$. The same inequality holds when the covariates are non-random; here we can view 
$\left(\eta(\bx_i)-\eta_0(\bx_i)\right)^2$; $i=1,\ldots,n$, as a set of independent realizations
from some independent stochastic process.
It follows that
\begin{align}
&\int_{\mathcal S^c}P\left(\left|\frac{1}{2\sigma^2}\times\frac{1}{n}\sum_{i=1}^n\left(\eta(\bx_i)-\eta_0(\bx_i)\right)^2
-\frac{1}{2\sigma^2}E_{\bX}\left(\eta(\bX)-\eta_0(\bX)\right)^2\right|>\frac{\kappa_1}{4}\right)d\pi(\theta)\notag\\
&\leq 2\int_{\mathcal G_n}\exp\left\{-\frac{nC\kappa^2_1\sigma^4}{\|\eta-\eta_0\|^2}\right\}d\pi(\theta)+\pi(\mathcal G^c_n)\notag\\
	&=2\int_{\|\eta\|\leq\exp(\left(\beta n\right)^{1/4})}\left[\int_{\exp(-2\left(\beta n\right)^{1/4})}^{\exp(2\left(\beta n\right)^{1/4})}
	\exp\left(-\frac{nC\kappa^2_1u^{-2}}{\|\eta-\eta_0\|^2}\right)\pi\left(u^{-1}\right)u^{-2}du\right]\pi\left(\|\eta\|\right)d\|\eta\|\notag\\
&\qquad+\pi(\mathcal G^c_n).
\label{eq:zero_prob3}
\end{align}

Replacing $\pi\left(u^{-1}\right)u^{-2}$ with $\tilde\pi_n(u)=c_3\sum_{r=1}^{M_n}\psi^{\zeta_{rn}}_{rn}\exp(-\psi_{rn} u^2)u^{2(\zeta_{rn}-1)}$, where
$M_n\leq\exp(\left(\beta n\right)^{1/4})$ is the number of mixture components, $c_3>0$, for $r=1,\ldots,M_n$, $\frac{1}{2}<\zeta_{rn}\leq c_4n^q$, for $0<q<1/4$ and $n\geq 1$, 
where $c_4>0$, and $0<\psi_1\leq\psi_{rn}<c_5<\infty$, for all $r$ and $n$, and using the same techniques as before, we obtain
\begin{align}
&\int_{\exp(-2\left(\beta n\right)^{1/4})}^{\exp(2\left(\beta n\right)^{1/4})}
	\exp\left(-\frac{nC\kappa^2_1u^{-2}}{\|\eta-\eta_0\|^2}\right)\pi\left(u^{-1}\right)u^{-2}du\notag\\
&\qquad\leq c_3\times\exp\left\{3\left(\beta n\right)^{1/4}+n^q\log c_5\right\}\times\exp\left\{-\frac{C_1\kappa_1\sqrt{\psi_1 n}}{\left(\|\eta\|+\|\eta_0\|\right)}\right\},	
\label{eq:new1}
\end{align}
for some constant $C_1>0$. 
Now, using the same techniques as before, we obtain
\begin{align}
	&\int_{\|\eta\|\leq\exp(\left(\beta n\right)^{1/4})} \exp\left\{-\frac{C_1\kappa_1\sqrt{\psi_1 n}}{\left(\|\eta\|+\|\eta_0\|\right)}\right\}
	\pi\left(\|\eta\|\right)d\|\eta\|\notag\\
	&\qquad=\int_{v\leq\|\eta_0\|+\exp(\left(\beta n\right)^{1/4})}\exp\left(-\frac{C_1\kappa_1\sqrt{\psi_1 n}}{v}\right)\pi\left(v-\|\eta_0\|\right)dv\notag\\
	&\qquad\leq c_3\sum_{r=1}^{M_n}c^{n^q}_5\int_{v\leq\|\eta_0\|+\exp\left(\left(\beta n\right)^{1/4}\right)}
	\exp\left\{-\left(\frac{C_1\kappa_1\sqrt{\psi_1 n}}{v}+\psi_1v-(\zeta_{rn}-1)\log v\right)\right\}dv\label{eq:new2}\\
	&\qquad\leq 2c_3\exp\left\{-\left(C_2\sqrt{\kappa_1}n^{1/4}-2\left(\beta n\right)^{1/4}-n^q\log c_5\right)\right\},\label{eq:new3}
\end{align}
with $\pi\left(v-\|\eta_0\|\right)$ replaced with the mixture as before.
Here $M_n\leq\exp(\left(\beta n\right)^{1/4})$, $c_3>0$, for $r=1,\ldots,M_n$, $0<\zeta_{rn}\leq c_4n^q$, for $0<q<1/4$ and $n\geq 1$, 
where $c_4>0$, and $0<\psi_1\leq\psi_{rn}<c_5<\infty$, for all $r$ and $n$.
Note that the negative of the exponent of the $r$-th term of (\ref{eq:new1}) is minimized for 
$\tilde v_{rn}=\frac{\zeta_{rn}-1+\sqrt{(\zeta_{rn}-1)^2+4\psi_1C_1\kappa_1\sqrt{\psi_1n}}}{2\psi_1}$, and for large $n$ it holds that
$\frac{\tilde C_1\sqrt{\kappa_1}n^{1/4}}{2\psi_1}\leq\tilde v_{rn}\leq \frac{\tilde C_2\sqrt{\kappa_1}n^{1/4}}{2\psi_1}$, for some positive constants
$\tilde C_1$ and $\tilde C_2$. Also, for large $n$, $\tilde v_{rn}\psi_1(1-c_2)\geq (\zeta_{rn}-1)\log\tilde v_{rn}$, for $0<c_2<1$. Hence (\ref{eq:new3})
follows from (\ref{eq:new2}) using $\frac{C_1\kappa_1\sqrt{\psi_1 n}}{\tilde v_{rn}}+\psi_1\tilde v_{rn}-(\zeta_{rn}-1)\log \tilde v_{rn}
\geq \psi_1\tilde v_{rn}-(\zeta_{rn}-1)\log \tilde v_{rn}\geq c_2\tilde v_{rn}\psi_1\geq C_2\sqrt{\kappa_1}n^{1/4}$, for some $C_2>0$.

Combining (\ref{eq:zero_prob3}), (\ref{eq:new1}) and (\ref{eq:new3}), we obtain
\begin{align}
&\int_{\mathcal S^c}P\left(\left|\frac{1}{2\sigma^2}\times\frac{1}{n}\sum_{i=1}^n\left(\eta(\bx_i)-\eta_0(\bx_i)\right)^2
	-\frac{1}{2\sigma^2}E_{\bX}\left(\eta(\bX)-\eta_0(\bX)\right)^2\right|>\frac{\kappa_1}{4}\right)d\pi(\theta)\notag\\
	&\qquad\leq C_1\exp\left\{-\left(C_2\sqrt{\kappa_1}n^{1/4}-5\left(\beta n\right)^{1/4}-2n^q\log c_5\right)\right\} + \pi\left(\mathcal G^c_n\right),
	\label{eq:new4}
\end{align}
where $C_1$ and $C_2$ are appropriate positive constants.
Since $\kappa_1$ is as large as desired, it follows that (\ref{eq:new3}) is summable.

For the fourth term, note that 
$$
Z_n=\frac{1}{n}\sum_{i=1}^n\left(\frac{y_i-\eta_0(\bx_i)}{\sigma_0}\right)\left(\eta(\bx_i)-\eta_0(\bx_i)\right) 
\sim N\left(0,\frac{1}{n^2}\sum_{i=1}^n\left(\eta(\bx_i)-\eta_0(\bx_i)\right)\right).$$
Then since
$$\sum_{i=1}^n\left(\eta(\bx_i)-\eta_0(\bx_i)\right)^2
\leq n\left(\underset{\bx\in\mathcal X}{\sup}~|\eta(\bx)-\eta_0(\bx)|\right)^2=n\|\eta-\eta_0\|^2,$$
\begin{align}
&P\left(\left|\frac{1}{\sigma^2}\times\frac{1}{n}\sum_{i=1}^n\left(y_i-\eta_0(\bx_i)\right)\left(\eta(\bx_i)-\eta_0(\bx_i)\right)\right|>\frac{\kappa_1}{4}\right)
= P\left(\left|Z_n\right|>\frac{\kappa_1\sigma^2}{4\sigma_0}\right)\notag\\
&\ \ \leq 2\exp\left(-\frac{Cn\kappa^2_1\sigma^4}{\sigma^2_0\|\eta-\eta_0\|^2}\right),
\label{eq:s5_9}
\end{align}
for some $C>0$.
Hence, in the same way as (\ref{eq:new4}), 
we obtain using (\ref{eq:s5_9}),
\begin{align}
&\int_{\mathcal S^c}P\left(\left|\frac{1}{\sigma^2}\times\frac{1}{n}\sum_{i=1}^n\left(y_i-\eta_0(\bx_i)\right)
\left(\eta(\bx_i)-\eta_0(\bx_i)\right)\right|>\frac{\kappa_1}{4}\right)d\pi(\theta)\notag\\
&\leq \int_{\mathcal G_n}2\exp\left(-\frac{Cn\kappa^2_1\sigma^4}{\sigma^2_0\|\eta-\eta_0\|^2}\right)d\pi(\theta)+2\pi\left(\mathcal G^c_n\right)\notag\\
&\leq C_1\exp\left\{-\left(C_2\sqrt{\kappa_1}n^{1/4}-5\left(\beta n\right)^{1/4}-2n^q\log c_5\right)\right\} + \pi\left(\mathcal G^c_n\right),
\label{eq:zero_prob5}
\end{align}
for relevant positive constants $C_1,C_2, c_5$.

Combining (\ref{eq:s5_6}), (\ref{eq:s5_7_int}), (\ref{eq:s5_8_int}), (\ref{eq:term1_bound3_2}), (\ref{eq:term1_bound3_3}), (\ref{eq:zero_prob3}), (\ref{eq:new4}), 
(\ref{eq:zero_prob5}),
and noting that $\sum_{n=1}^{\infty}\pi\left(\mathcal G^c_n\right)<\sum_{n=1}^{\infty}\alpha\exp\left(-\beta n\right)<\infty$, we obtain
\begin{equation*}
\int_{\mathcal S^c}P\left(\left|\frac{1}{n}\log R_n(\theta)+h(\theta)\right|>\kappa_1\right)d\pi(\theta)<\infty.
\end{equation*}

\subsection{Verification of (S7)}
\label{subsec:S7}
For any set $A$ such that $\pi(A)>0$, $\mathcal G_n\cap A\uparrow A$. It follows from this and continuity of $h$ that $h\left(\mathcal G_n\cap A\right)\downarrow h\left(A\right)$ as
$n\rightarrow\infty$, so that (S7) holds.


\section{Verification of Shalizi's conditions for Gaussian process regression with double exponential error distribution}
\label{sec:de}



\subsection{Verification of (S1)}
\label{subsec:S1_de}
In this case,
\begin{equation}
\frac{1}{n}\log R_n(\theta)=\log\left(\frac{\sigma_0}{\sigma}\right)+\frac{1}{\sigma_0}\times\frac{1}{n}\sum_{i=1}^n\left|y_i-\eta_0(\bx_i)\right)|
-\frac{1}{\sigma}\times\frac{1}{n}\sum_{i=1}^n\left|y_i-\eta(\bx_i)\right|.
\label{eq:R_de}
\end{equation}
%
As before, note that
\begin{align}
&\Bigg|\frac{1}{n}\sum_{i=1}^n\left|y_{1i}-\eta_0(\bx_i)\right|-\frac{1}{n}\sum_{i=1}^n\left|y_{2i}-\eta_0(\bx_i)\right|\Bigg |\notag\\
&\leq\frac{1}{n}\sum_{i=1}^n\Bigg |\left|y_{1i}-\eta_0(\bx_i)\right|-\left|y_{2i}-\eta_0(\bx_i)\right|\Bigg |\notag\\
&\leq\frac{1}{n}\sum_{i=1}^n|y_{1i}-y_{2i}|\notag\\
&\leq n^{-\frac{1}{2}}\sqrt{\sum_{i=1}^n(y_{1i}-y_{2i})^2}\notag\\
&=n^{-\frac{1}{2}}\|\by_{1n}-\by_{2n}\|,\notag
\end{align}
from which Lipschitz continuity follows.
Similarly, 
\begin{align}
&\Bigg|\frac{1}{n}\sum_{i=1}^n\left|y_{1i}-\eta_1(\bx_i)\right|-\frac{1}{n}\sum_{i=1}^n\left|y_{2i}-\eta_2(\bx_i)\right|\Bigg |\notag\\ 
&\leq\frac{1}{n}\sum_{i=1}^n\left|y_{1i}-\eta_1(\bx_i)-y_{2i}+\eta_2(\bx_i)\right|\notag\\
&\leq\frac{1}{n}\sum_{i=1}^n\left[\left|y_{1i}-y_{2i}\right|+\left|\eta_1(\bx_i)-\eta_2(\bx_i)\right|\right]\notag\\
&\leq n^{-\frac{1}{2}}\|\by_1-\by_2\|+\|\eta_1-\eta_2\|,
\label{eq:lipcont_de}
\end{align}
which implies continuity of $\frac{1}{n}\sum_{i=1}^n\left|y_i-\eta(\bx_i)\right|$ with respect to $\by$ and $\eta$.
In other words, (\ref{eq:R2}) is continuous and hence measurable, as before. Measurability, when the covariates are considered random, also follows as before,
using measurability of $\eta_0(\bX)$ as assumed in (A4).

\subsection{Verification of (S2) and proof of Lemma \ref{lemma:lemma1} for double-exponential errors}
\label{subsec:S2_de}
Now note that if $\epsilon_i=y_i-\eta_0(\bx_i)$ has the double exponential density of the form 
\begin{equation*}
f(\epsilon)=\frac{1}{2\sigma}\exp\left(-\frac{|\epsilon|}{\sigma}\right);~\epsilon\in\mathbb R.
\end{equation*}
with $\sigma$ replaced with $\sigma_0$, then
\begin{align}
&E_{\theta_0}\left|y_i-\eta_0(\bx_i)\right|=\sigma_0;\label{eq:de_exp1}\\ 
&E_{\theta_0}\left|y_i-\eta(\bx_i)\right|=E_{\theta_0}\left|(y_i-\eta_0(\bx_i))+(\eta_0(\bx_i)-\eta(\bx_i))\right|\notag\\
&\qquad\qquad\qquad\quad=|\eta_0(\bx_i)-\eta(\bx_i)|+\sigma_0\exp\left(-\frac{|\eta_0(\bx_i)-\eta(\bx_i)|}{\sigma_0}\right).
\label{eq:de_exp2}
\end{align}
It follows from (\ref{eq:de_exp1}), (\ref{eq:de_exp2}) and (A3), that
\begin{align}
&\frac{1}{n}\sum_{i=1}^nE_{\theta_0}\left|y_i-\eta_0(\bx_i)\right|=\sigma_0;\label{eq:de_exp3}\\
&\frac{1}{n}\sum_{i=1}^nE_{\theta_0}\left|y_i-\eta(\bx_i)\right|=
\frac{1}{n}\sum_{i=1}^n\left[|\eta(\bx_i)-\eta_0(\bx_i)|+\sigma_0\exp\left(-\frac{|\eta(\bx_i)-\eta_0(\bx_i)|}{\sigma_0}\right)\right]\notag\\
&\ \ \rightarrow E_{\bX}\left|\eta(\bX)-\eta_0(\bX)\right|+\sigma_0E_{\bX}\left[\exp\left(-\frac{|\eta(\bX)-\eta_0(\bX)|}{\sigma_0}\right)\right],
~\mbox{as}~n\rightarrow\infty.\label{eq:de_exp4}
\end{align}
Using (\ref{eq:de_exp3}) and (\ref{eq:de_exp4}) we see that as $n\rightarrow\infty$,
\begin{align}
&\frac{1}{n}E_{\theta_0}\left[\log R_n(\theta)\right]=\log\left(\frac{\sigma_0}{\sigma}\right)
+\frac{1}{\sigma_0}\times\frac{1}{n}\sum_{i=1}^nE_{\theta_0}\left|y_i-\eta_0(\bx_i)\right|
-\frac{1}{\sigma}\times\frac{1}{n}\sum_{i=1}^nE_{\theta_0}\left|y_i-\eta(\bx_i)\right|\notag\\
&\rightarrow\log\left(\frac{\sigma_0}{\sigma}\right)+1-\frac{1}{\sigma}E_{\bX}\left|\eta(\bX)-\eta_0(\bX)\right|
-\frac{\sigma_0}{\sigma}E_{\bX}\left[\exp\left(-\frac{|\eta(\bX)-\eta_0(\bX)|}{\sigma_0}\right)\right],\notag\\
&\ \ =-h(\theta),
\label{eq:R_de2}
\end{align}
where
\begin{equation*}
h(\theta)=\log\left(\frac{\sigma}{\sigma_0}\right)-1+\frac{1}{\sigma}E_{\bX}\left|\eta(\bX)-\eta_0(\bX)\right|
+\frac{\sigma_0}{\sigma}E_{\bX}\left[\exp\left(-\frac{|\eta(\bX)-\eta_0(\bX)|}{\sigma_0}\right)\right].
\end{equation*}
As in the case of Gaussian errors, the results remain the same if the covariates are assumed to be random.

\subsection{Verification of (S3) and proof of Theorem \ref{theorem:theorem1} for double exponential errors}
\label{subsec:S3_de}
We now show that for all $\theta\in\Theta$, $\underset{n\rightarrow\infty}{\lim}\frac{1}{n}\log R_n(\theta)=-h(\theta)$, almost surely.
First note that 
\begin{align}
&\left|\frac{1}{n}R_n(\theta)+h(\theta)\right|\leq\left|\frac{1}{n}\sum_{i=1}^n\frac{|y_i-\eta_0(\bx_i)|}{\sigma_0}-1\right|\notag\\
&\quad+\left|\frac{1}{n}\sum_{i=1}^n\frac{|y_i-\eta(\bx_i)|}{\sigma}-\frac{1}{\sigma}E_{\bX}\left|\eta(\bX)-\eta_0(\bX)\right|-
\frac{\sigma_0}{\sigma}E_{\bX}\left[\exp\left(-\frac{|\eta(\bX)-\eta_0(\bX)|}{\sigma_0}\right)\right]\right|.
\label{eq:h_de2}
\end{align}
Since $\frac{|y_i-\eta_0(\bx_i)|}{\sigma_0}$ has the exponential distribution with mean one, the term 
$\left|\frac{1}{n}\sum_{i=1}^n\frac{|y_i-\eta_0(\bx_i)|}{\sigma_0}-1\right|\rightarrow 0$ almost surely as $n\rightarrow\infty$ by the
strong law of large numbers.
That the term (\ref{eq:h_de2}) also tends to zero almost surely as $n\rightarrow\infty$ can be shown using the Borel-Cantelli lemma, using the inequality
(\ref{eq:bound3}), and replacing $\kappa_1$ in that inequality with any $\delta_1>0$.
In other words, it holds that for all $\theta\in\Theta$, $\underset{n\rightarrow\infty}{\lim}\frac{1}{n}\log R_n(\theta)=-h(\theta)$, almost surely.
Also, it follows from (\ref{eq:R_de}), (\ref{eq:lipcont_de}), (\ref{eq:h_de}), Lipschitz continuity of $x\mapsto\exp(-|x|)$, boundedness of the first derivative with respect to $\sigma$, 
that $\frac{1}{n}\log R_n(\theta)+h(\theta)$ is Lipschitz on $\theta\in\mathcal G_n\setminus I=\mathcal G_n$, which is compact. 
As a result, it follows that $\frac{1}{n}\log R_n(\theta)+h(\theta)$ is stochastically equicontinuous in $\mathcal G\in\left\{\mathcal G_1,\mathcal G_2,\ldots,\right\}$. Hence, the convergence  
$\underset{n\rightarrow\infty}{\lim}\frac{1}{n}\log R_n(\theta)=-h(\theta)$ occurs uniformly for $\theta\in\mathcal G$, almost surely.

\subsection{Verification of (S4)}
\label{subsec:S4_de}
Note that $h(\theta)\leq\log\left(\frac{\sigma}{\sigma_0}\right)-1+\frac{\|\eta-\eta_0\|+\sigma_0}{\sigma}$. Now $0<\|\eta-\eta_0\|<\infty$ and 
$0<\sigma<\infty$ with prior probability one. Consequently, it follows that
$h(\theta)<\infty$ with probability one, so that $I=\emptyset$ and hence, $\mathcal G_n\setminus I=\mathcal G_n$.

\subsection{Verification of (S5)}
\label{subsec:S5_de}

Verification of (S5) (1) and (S5) (2) remains the same as for Gaussian noise. (S5) (3) follows in the same way as for Gaussian noise is we can show that
$h(\theta)$ is continuous in $\theta$.
To see that $h(\theta)$ is continuous in $\theta$, again assume that $\eta_j\rightarrow\tilde\eta$ as $j\rightarrow\infty$ in the sense that $\|\eta_j-\tilde\eta\|\rightarrow 0$
as $j\rightarrow\infty$. Then $\left|E_{\bX}\left|\eta_j(\bX)-\eta_0(\bX)\right|-E_{\bX}\left|\tilde\eta(\bX)-\eta_0(\bX)\right|\right|\leq E_{\bX}\left|\eta_j(\bX)-\tilde\eta(\bX)\right|
\leq\|\eta_j-\tilde\eta\|\rightarrow 0$ as $j\rightarrow\infty$. Also, 
\begin{align}
&\left|E_{\bX}\left[\exp\left(-\frac{|\eta_j(\bX)-\eta_0(\bX)|}{\sigma_0}\right)\right]-E_{\bX}\left[\exp\left(-\frac{|\tilde\eta(\bX)-\eta_0(\bX)|}{\sigma_0}\right)\right]\right|\notag\\ 
&\leq E_{\bX}\left[\exp\left(-\left|\tilde\eta(\bX)-\eta_0(\bX)\right|\right)\times
\left|\exp\left(-\frac{\left(|\eta_j(\bX)-\eta_0(\bX)|-|\tilde\eta(\bX)-\eta_0(\bX)|\right)}{\sigma_0}\right)-1\right|\right]\notag\\
&\leq E_{\bX}\left[\exp\left(-\left|\tilde\eta(\bX)-\eta_0(\bX)\right|\right)\times
\left|\exp\left(\frac{|\eta_j(\bX)-\tilde\eta(\bX)|}{\sigma_0}\right)-1\right|\right]\notag\\
&\leq \left|\exp\left(\frac{\|\eta_j-\tilde\eta\|}{\sigma_0}\right)-1\right|\times E_{\bX}\left[\exp\left(-\left|\tilde\eta(\bX)-\eta_0(\bX)\right|\right)\right]\notag\\
&\rightarrow 0,~\mbox{as}~j\rightarrow\infty.\notag
\end{align}
Continuity of $h(\theta)$ hence follows easily.

\subsection{Verification of (S6) and proof of Theorem \ref{theorem:theorem3} for double exponential errors}
\label{subsec:S6_de}
It follows from (\ref{eq:h_de2}) that for all $\theta\in\Theta$, for $\kappa_1=\kappa-h(\Theta)$, we have
\begin{align}
&P\left(\left|\frac{1}{n}\log R_n(\theta)+h(\theta)\right|>\kappa_1\right)
\leq P\left(\left|\frac{1}{\sigma_0}\times\frac{1}{n}\sum_{i=1}^n\left|y_i-\eta_0(\bx_i)\right|-1\right|>\frac{\kappa_1}{2}\right)\notag\\
&\qquad+P\left(\left|\frac{1}{\sigma}\times\frac{1}{n}\sum_{i=1}^n\left|y_i-\eta(\bx_i)\right|
-\frac{1}{\sigma}E_{\bX}\left|\eta(\bX)-\eta_0(\bX)\right|\right.\right.\notag\\
&\qquad\qquad\qquad\left.\left.-\frac{\sigma_0}{\sigma}E_{\bX}\left(\exp\left\{-\frac{\left|\eta(\bX)-\eta_0(\bX)\right|}{\sigma_0}\right\}\right)\right|
>\frac{\kappa_1}{2}\right).
\label{eq:unif_de1}
\end{align}

Since $\frac{|y_i-\eta_0(\bx_i)|}{\sigma_0}$ are exponential random variables with expectation one, it follows that $\frac{|y_i-\eta_0(\bx_i)|}{\sigma_0}-1$ are zero-mean, 
independent sub-exponential random variables with some parameter $s>0$. Hence, by Bernstein's inequality (\ctn{Uspensky37}, \ctn{Bennett62}, \ctn{Massart03}),
\begin{align}
&P\left(\left|\frac{1}{\sigma_0}\times\frac{1}{n}\sum_{i=1}^n\left|y_i-\eta_0(\bx_i)\right|-1\right|>\frac{\kappa_1}{2}\right)
\leq 2\exp\left(-\frac{n}{2}\min\left\{\frac{\kappa^2_1}{4s^2},\frac{\kappa_1}{2s}\right\}\right).\notag
\end{align}
Hence,
\begin{equation}
\int_{\mathcal S^c}P\left(\left|\frac{1}{\sigma_0}\times\frac{1}{n}\sum_{i=1}^n\left|y_i-\eta_0(\bx_i)\right|-1\right|>\frac{\kappa_1}{2}\right)
\leq 2\exp\left(-\frac{n}{2}\min\left\{\frac{\kappa^2_1}{4s^2},\frac{\kappa_1}{2s}\right\}\right).
\label{eq:summable1}
\end{equation}

Let $\bar\varphi=E_{\bX}\left|\eta(\bX)-\eta_0(\bX)\right|+\sigma_0E_{\bX}\left(\exp\left\{-\frac{\left|\eta(\bX)-\eta_0(\bX)\right|}{\sigma_0}\right\}\right)$. 
Also, letting
$\varphi(\bx)=\left|\eta(\bx)-\eta_0(\bx)\right|+\sigma_0\left(\exp\left\{-\frac{\left|\eta(\bx)-\eta_0(\bx)\right|}{\sigma_0}\right\}\right)$, note that
\begin{equation}
\frac{1}{n}\sum_{i=1}^n\varphi(\bx_i)\rightarrow\bar\varphi,~\mbox{as}~n\rightarrow\infty.
\label{eq:varphi1}
\end{equation}
With this, the second term of (\ref{eq:unif_de1}) can be bounded as follows:
\begin{align}
&P\left(\left|\frac{1}{\sigma}\times\frac{1}{n}\sum_{i=1}^n\left|y_i-\eta(\bx_i)\right|-\frac{\bar\varphi}{\sigma}\right|>\frac{\kappa_1}{2}\right)\notag\\
&\ \ = P\left(\sigma^{-1}\left|\frac{1}{n}\sum_{i=1}^n\left\{\left|y_i-\eta(\bx_i)\right|-\varphi(\bx_i)\right\}+\frac{1}{n}\sum_{i=1}^n\varphi(\bx_i)-\bar\varphi\right|
>\frac{\kappa_1}{2}\right)\notag\\
&\ \ \leq P\left(\sigma^{-1}\left|\frac{1}{n}\sum_{i=1}^n\left\{\left|y_i-\eta(\bx_i)\right|-\varphi(\bx_i)\right\}\right|>\frac{\kappa_1}{4}\right)
+P\left(\sigma^{-1}\left|\frac{1}{n}\sum_{i=1}^n\varphi(\bx_i)-\bar\varphi\right|>\frac{\kappa_1}{4}\right).
\label{eq:bound1}
\end{align}
%
In the case of random or non-random covariates $\bX$, again by Hoeffding's inequality,
\begin{equation}
P\left(\sigma^{-1}\left|\frac{1}{n}\sum_{i=1}^n\varphi(\bx_i)-\bar\varphi\right|>\frac{\kappa_1}{4}\right)
\leq\exp\left\{-\frac{nC\kappa^2_1\sigma^2}{\left(\|\eta\|+c_0\right)^2}\right\},
\label{eq:hoeff1}
\end{equation}
where $C>0$ is a constant. Note that $\left(\|\eta\|+c_0\right)$, with $c_0=\|\eta_0\|+\sigma_0$, is an upper bound of $|\varphi(\cdot)|$. 
Again, such an upper bound exists since $\mathcal X$ is compact and $\eta(\cdot)$ is continuous on $\mathcal X$.
Application of the same method as proving 
(\ref{eq:term1_bound3_3}) and (\ref{eq:new4}) yields
\begin{align}
&\int_{\mathcal S^c}P\left(\sigma^{-1}\left|\frac{1}{n}\sum_{i=1}^n\varphi(\bx_i)-\bar\varphi\right|>\frac{\kappa_1}{4}\right)\pi(\theta)d\theta\notag\\
&\qquad\leq C_1\exp\left\{-\left(C_2\sqrt{\kappa_1}n^{1/4}-5\left(\beta n\right)^{1/4}-2n^q\log c_5\right)\right\} + \pi\left(\mathcal G^c_n\right),
\label{eq:summable2}
\end{align}
where as before $\kappa_1$ is large enough to make the exponent of (\ref{eq:summable2}) negative.

For the first term of (\ref{eq:bound1}), let us first prove that $|y_i-\eta(\bx_i)|-\varphi(\bx_i)$ are sub-exponential random variables. Then we can apply Bernstein's
inequality to directly bound the term. 
We need to show that $E_{\theta_0}\left[\exp\left\{t\left(|y_i-\eta(\bx_i)|\right)-\varphi(\bx_i)\right\}\right]\leq\exp\left(\frac{t^2s^2}{2}\right)$
for $|t|\leq s^{-1}$, for some $s>0$.

\subsection{Case 1: $t\geq 0$, $\eta(\bx_i)-\eta_0(\bx_i)>0$}
\label{subsec:case1}
Direct calculation shows that
\begin{align}
&E_{\theta_0}\left[\exp\left\{t\left(|y_i-\eta(\bx_i)|\right)-\varphi(\bx_i)\right\}\right]\notag\\
&=\exp\left(-t\varphi(\bx_i)\right)\times
\frac{\exp\left\{(\eta(\bx_i)-\eta_0(\bx_i))t\right\}-\exp\left(\frac{\eta_0(\bx_i)-\eta(\bx_i)}{\sigma_0}\right)}{1-\sigma^2_0t^2}\notag\\
&\leq\frac{\exp\left\{t\left(\varphi(\bx_i)+\eta(\bx_i)-\eta_0(\bx_i)\right)\right\}}{1-\sigma^2_0t^2}\notag\\
&\leq \frac{\exp\left\{t\left(2\|\eta-\eta_0\|+\sigma_0\right)\right\}}{1-\sigma^2_0t^2}.
\label{eq:subexp1}
\end{align}
To show that (\ref{eq:subexp1}) is bounded above by $\exp(t^2s^2/2)$, we need to show that 
\begin{equation}
f(t)=\frac{t^2s^2}{2}-2(\|\eta-\eta_0\|+\sigma_0)t+\log(1-\sigma^2_0t^2)\geq 0.
\label{eq:subexp2}
\end{equation}
For $t>0$, it is sufficient to show that 
\begin{equation}
\frac{ts^2}{2}\geq 2(\|\eta-\eta_0\|+\sigma_0)-\frac{\log(1-\sigma^2_0t^2)}{t}.
\label{eq:subexp3}
\end{equation}
Now, $-\frac{\log(1-\sigma^2_0t^2)}{t}\rightarrow 0$, as $t\rightarrow 0$. Hence, for any $\epsilon>0$, there exists $\delta(\epsilon)>0$ such that $t\leq\delta(\epsilon)$
implies $-\frac{\log(1-\sigma^2_0t^2)}{t}<\epsilon$. Let $s\geq\frac{C_1\|\eta-\eta_0\|+C_2}{\delta(\epsilon)}$, where $C_1>0$ and $C_2>0$ are sufficiently large quantities. 
Hence, if $\delta(\epsilon)^2\leq t\leq\delta(\epsilon)$, then (\ref{eq:subexp3}), and hence (\ref{eq:subexp2}), is satisfied. Now, $f(t)$ given by (\ref{eq:subexp2}) is continuous in $t$
and $f(0)=0$. Hence, (\ref{eq:subexp2}) holds even for $0\leq t\leq\delta(\epsilon)^2$. In other words,
\begin{equation}
E_{\theta_0}\left[\exp\left\{t\left(|y_i-\eta(\bx_i)|\right)-\varphi\right\}\right]\leq\exp\left(\frac{t^2s^2}{2}\right),~\mbox{for}~0\leq t\leq s^{-1}
\leq\frac{\delta(\epsilon)}{C_1\|\eta-\eta_0\|+C_2}\leq\delta(\epsilon).
\label{eq:subexp4}
\end{equation}

\subsection{Case 2: $t\geq 0$, $\eta(\bx_i)-\eta_0(\bx_i)<0$}
\label{subsec:case2}
In this case,
\begin{align}
&E_{\theta_0}\left[\exp\left\{t\left(|y_i-\eta(\bx_i)|\right)-\varphi(\bx_i)\right\}\right]\notag\\
&=\exp\left(-t\varphi(\bx_i)\right)\times
\frac{\exp\left\{(\eta_0(\bx_i)-\eta(\bx_i))t\right\}+\sigma_0t\exp\left(\frac{\eta(\bx_i)-\eta_0(\bx_i)}{\sigma_0}\right)}{1-\sigma^2_0t^2}\notag\\
&\leq \exp\left(t\varphi(\bx_i)\right)\times
\frac{2\exp\left\{(\eta_0(\bx_i)-\eta(\bx_i))t\right\}}{1-\sigma^2_0t^2}\notag\\
&\leq\frac{\exp\left\{t\left(\varphi(\bx_i)+(\eta_0(\bx_i)-\eta(\bx_i))\right)\right\}}{\frac{1-\sigma^2_0t^2}{2}}\notag\\
&\leq \frac{\exp\left\{t\left(2\|\eta-\eta_0\|+\sigma_0\right)\right\}}{\frac{1-\sigma^2_0t^2}{2}}.\notag
\end{align}
As in Section \ref{subsec:case1} it can be seen that (\ref{eq:subexp4}) holds.

\subsection{Case 3: $t\leq 0$, $\eta(\bx_i)-\eta_0(\bx_i)>0$}
\label{subsec:case3}
Here
\begin{align}
&E_{\theta_0}\left[\exp\left\{t\left(|y_i-\eta(\bx_i)|\right)-\varphi(\bx_i)\right\}\right]\notag\\
&=\exp\left(-t\varphi(\bx_i)\right)\times
\frac{\exp\left\{(\eta(\bx_i)-\eta_0(\bx_i))t\right\}-\sigma_0|t|\exp\left(\frac{\eta_0(\bx_i)-\eta(\bx_i)}{\sigma_0}\right)}{1-\sigma^2_0t^2}\notag\\
&\leq \exp\left(-t\varphi(\bx_i)\right)\times
\frac{1}{1-\sigma^2_0t^2}\notag\\
&\leq \frac{\exp\left\{-t\left(\|\eta-\eta_0\|+\sigma_0\right)\right\}}{1-\sigma^2_0t^2}.\notag
\end{align}
Here we need to have
$|t|\left[\frac{|t|s^2}{2}-\left(\|\eta-\eta_0\|+\sigma_0\right)+\frac{\log(1-\sigma^2_0t^2)}{|t|}\right]>0$. In the same way as before it follows that
\begin{equation}
E_{\theta_0}\left[\exp\left\{t\left(|y_i-\eta(\bx_i)|\right)-\varphi\right\}\right]\leq\exp\left(\frac{t^2s^2}{2}\right),~\mbox{for}~0\leq |t|\leq s^{-1}
\leq\frac{\delta(\epsilon)}{C_1\|\eta-\eta_0\|+C_2}\leq\delta(\epsilon).\notag
\label{eq:subexp7}
\end{equation}

\subsection{Case 4: $t\leq 0$, $\eta(\bx_i)-\eta_0(\bx_i)<0$}
\label{subsec:case4}
In this case,
\begin{align}
&E_{\theta_0}\left[\exp\left\{t\left(|y_i-\eta(\bx_i)|\right)-\varphi(\bx_i)\right\}\right]\notag\\
&=\exp\left(-t\varphi(\bx_i)\right)\times
\frac{\exp\left\{(\eta_0(\bx_i)-\eta(\bx_i))t\right\}-\sigma_0|t|\exp\left(\frac{\eta(\bx_i)-\eta_0(\bx_i)}{\sigma_0}\right)}{1-\sigma^2_0t^2}\notag\\
&\leq \exp\left(-t\varphi(\bx_i)\right)\times
\frac{1}{1-\sigma^2_0t^2}\notag\\
&\leq \frac{\exp\left\{-t\left(\|\eta-\eta_0\|+\sigma_0\right)\right\}}{1-\sigma^2_0t^2}.\notag
\end{align}
Hence, (\ref{eq:subexp7}) holds.

Hence, for $i=1,\ldots,n$, $|y_i-\eta(\bx_i)|-E\left(|y_i-\eta(\bx_i)|\right)$ are 
zero-mean, independent sub-exponential random variables with parameter $s$. In particular, we can set $s=\frac{C_1\|\eta-\eta_0\|+C_2}{\delta(\epsilon)}$. Hence, by Bernstein's inequality,
\begin{align}
&P\left(\sigma^{-1}\left|\frac{1}{n}\sum_{i=1}^n\left\{\left|y_i-\eta(\bx_i)\right|-\varphi(\bx_i)\right\}\right|>\frac{\kappa_1}{4}\right)\notag\\
&\leq 2\max\left\{P\left(\frac{\sigma^{-1}}{n}\sum_{i=1}^n\left\{\left|y_i-\eta(\bx_i)\right|-\varphi(\bx_i)\right\}>\frac{\kappa_1}{4}\right),
P\left(\frac{\sigma^{-1}}{n}\sum_{i=1}^n\left\{\left|y_i-\eta(\bx_i)\right|-\varphi(\bx_i)\right\}< -\frac{\kappa_1}{4}\right)\right\}\notag\\
&\leq 2\exp\left(-\frac{n}{2}\min\left\{\frac{\kappa^2_1\sigma^2}{16s^2},\frac{\kappa_1\sigma}{4s}\right\}\right)\notag\\
&=2\exp\left(-\frac{n}{2}\min\left\{\frac{\kappa^2_1\delta(\epsilon)^2\sigma^2}{16(C_1\|\eta-\eta_0\|+C_2)^2},\frac{\kappa_1\delta(\epsilon)\sigma}{4(C_1\|\eta-\eta_0\|+C_2)}\right\}\right).
\label{eq:bound3}
\end{align}
Hence,
\begin{align}
&\int_{\mathcal S^c}P\left(\sigma^{-1}\left|\frac{1}{n}\sum_{i=1}^n\left\{\left|y_i-\eta(\bx_i)\right|-\varphi(\bx_i)\right\}\right|>\frac{\kappa_1}{4}\right)
	d\pi(\theta)\notag\\
	&\qquad\leq\int_{\mathcal G_n}2\exp\left(-\frac{n\kappa^2_1\delta(\epsilon)^2\sigma^2}{32(C_1\|\eta-\eta_0\|+C_2)^2}\right)d\pi(\theta)
	+\int_{\mathcal G_n}2\exp\left(-\frac{n\kappa_1\delta(\epsilon)\sigma}{8(C_1\|\eta-\eta_0\|+C_2)}\right)\label{eq:new1_1}\\
	&\qquad\qquad+2\pi\left(\mathcal G^c_n\right).\label{eq:new1_2}
\end{align}	

Applying the same techniques as proving (\ref{eq:new4}) we obtain
\begin{align}
&\int_{\mathcal G_n}2\exp\left(-\frac{n\kappa^2_1\delta(\epsilon)^2\sigma^2}{32(C_1\|\eta-\eta_0\|+C_2)^2}\right)d\pi(\theta)
\leq C_1\exp\left\{-\left(C_2\sqrt{\kappa_1}n^{1/4}-5\left(\beta n\right)^{1/4}-2n^q\log c_5\right)\right\},
\label{eq:new1_3}
\end{align}
for appropriate positive constants $C_1,C_2,c_5$.

For the second integral of (\ref{eq:new1_1}), observe that for appropriate positive constant $c_0$ and $C$,
\begin{align}
&\int_{\mathcal G_n}2\exp\left(-\frac{n\kappa_1\delta(\epsilon)\sigma}{8(C_1\|\eta-\eta_0\|+C_2)}\right)\notag\\
	&\leq 2\int_{\|\eta\|\leq\exp\left(\left(\beta n\right)^{1/4}\right)}
	\left[\int_{\exp\left(-2\left(\beta n\right)^{1/4}\right)}^{\exp\left(-2\left(\beta n\right)^{1/4}\right)}
	\exp\left(-\frac{C\kappa_1nu^{-1}}{\|\eta\|+c_0}\right)2\pi(u^{-2})u^{-3}du\right]\pi\left(\|\eta\|\right)d\|\eta\|
	\label{eq:new1_4}
\end{align}
Replacing $2\pi(u^{-2})u^{-3}$ with the mixture form as before with $0<\zeta_{rn}<c_5n^q$, where $0<q<1/4$, and the rest remaing the same as before, we obtain
\begin{equation}
	\int_{\exp\left(-2\left(\beta n\right)^{1/4}\right)}^{\exp\left(-2\left(\beta n\right)^{1/4}\right)}
	\exp\left(-\frac{C\kappa_1nu^{-1}}{\|\eta\|+c_0}\right)2\pi(u^{-2})u^{-3}du
	\leq\exp\left(-\frac{C_1\kappa_1\sqrt{n}}{\sqrt{\|\eta\|+c_0}}\right),
	\label{eq:new1_5}
\end{equation}
for some appropriate positive constant $C_1$.

Now we obtain
\begin{align}
&\int_{\|\eta\|\leq\exp\left(\left(\beta n\right)^{1/4}\right)}\exp\left(-\frac{C_1\kappa_1\sqrt{n}}{\sqrt{\|\eta\|+c_0}}\right)\pi\left(\|\eta\|\right)d\|\eta\|\notag\\
	&=\int_{0\leq v\leq\sqrt{c_0+\exp\left(\left(\beta n\right)^{1/4}\right)}}\exp\left(-\frac{C_2\kappa_1\sqrt{n}}{v}\right)\pi\left(v^2-c_0\right)2vdv\notag\\
	&\leq \tilde C_1\exp\left\{-\left(\tilde C_2\sqrt{\kappa_1}n^{1/4}-\frac{9}{2}\left(\beta n\right)^{1/4}-2n^q\log c_5\right)\right\},\label{eq:new1_6}
\end{align}
for appropriate positive constants $\tilde C_1$ and $\tilde C_2$. From (\ref{eq:new1_5}) and (\ref{eq:new1_6}) it follows that (\ref{eq:new1_6}) 
is an upper bound for (\ref{eq:new1_4}).
Combining this with (\ref{eq:bound3}), (\ref{eq:new1_1}), (\ref{eq:new1_2}) and (\ref{eq:new1_3}), we obtain
\begin{align}
&P\left(\sigma^{-1}\left|\frac{1}{n}\sum_{i=1}^n\left\{\left|y_i-\eta(\bx_i)\right|-\varphi(\bx_i)\right\}\right|>\frac{\kappa_1}{4}\right)\notag\\
	&\leq C_1\exp\left\{-\left(C_2\sqrt{\kappa_1}n^{1/4}-5\left(\beta n\right)^{1/4}-2n^q\log c_5\right)\right\}\notag\\
	&\qquad+\tilde C_1\exp\left\{-\left(\tilde C_2\sqrt{\kappa_1}n^{1/4}-\frac{9}{2}\left(\beta n\right)^{1/4}-2n^q\log c_5\right)\right\}+2\pi\left(\mathcal G^c_n\right).
\label{eq:summable3}
\end{align}

Gathering (\ref{eq:summable1}), (\ref{eq:summable2}) and (\ref{eq:summable3}) we see that
\begin{equation}
\sum_{n=1}^{\infty}\int_{\mathcal S^c}P\left(\left|\frac{1}{n}\log R_n(\theta)+h(\theta)\right|>\delta\right)\pi(\theta)d\theta<\infty.
\label{eq:unif_de2}
\end{equation}

\subsection{Verification of (S7)}
\label{subsec:S7_de}
Verification of (S7) is exactly the same as for Gaussian errors.

\section{Verification of the assumptions of Shalizi for the general stochastic process model}
\label{sec:nongp_verification}

Note that
\begin{align}
f_{\theta}(\by_n)&=\frac{1}{\sigma^n}\prod_{i=1}^n\phi(y_i-\eta(\bx_i));\label{eq:nongp_like1}\\
f_{\theta_0}(\by_n)&=\frac{1}{\sigma^n_0}\prod_{i=1}^n\phi(y_i-\eta_0(\bx_i)).\label{eq:nongp_true_like1}
\end{align}

\subsection{Verification of (S1)}
\label{subsec:nongp_S1}
The equations (\ref{eq:nongp_like1}) and (\ref{eq:nongp_true_like1}) yield, in our case,
\begin{equation}
\frac{1}{n}\log R_n(\theta)=\log\left(\frac{\sigma_0}{\sigma}\right)+\frac{1}{n}\sum_{i=1}^n\log\phi\left(\frac{y_i-\eta_0(\bx_i)}{\sigma_0}\right)
-\frac{1}{n}\sum_{i=1}^n\log\phi\left(\frac{y_i-\eta(\bx_i)}{\sigma}\right).
\label{eq:nongp_R1}
\end{equation}
We show that the right hand side of (\ref{eq:nongp_R1}), which we denote as $f(\by_n,\theta)$, is continuous in $(\by_n,\theta)$, which is sufficient to confirm measurability of $R_n(\theta)$.
Let $\|(\by_n,\theta)\|=\|\by_n\|+\|\theta\|$, where $\|\by_n\|$ is the Euclidean norm and $\|\theta\|=\|\eta\|+|\sigma|$, with
$\|\eta\|=\underset{\bx\in\mathcal X}{\sup}~|\eta(\bx)|$. Since $\mathcal X$ is compact and $\eta$ is almost surely continuous, it follows that $\|\eta\|<\infty$ almost surely.

Consider $\by_{1n}=(y_{11},y_{12},\ldots,y_{1n})^T$, $\by_{2n}=(y_{21},y_{22},\ldots,y_{2n})^T$, $\theta_1$ and $\theta_2$. 
Using the Lipschitz condition of (A7), we obtain
\begin{align}
&\Bigg|\frac{1}{n}\sum_{i=1}^n\log\phi\left(\frac{y_{1i}-\eta_0(\bx_i)}{\sigma_0}\right)-\frac{1}{n}\sum_{i=1}^n\log\phi\left(\frac{y_{2i}-\eta_0(\bx_i)}{\sigma_0}\right)\Bigg |\label{eq:cont1}\\
&\leq\frac{L}{n\sigma_0}\sum_{i=1}^n|y_{1i}-y_{2i}|\leq \frac{L}{n\sigma_0}\|\by_{1n}-\by_{2n}\|.
\label{eq:nongp_cont2}
\end{align}
Hence, 
the term $\frac{1}{n}\sum_{i=1}^n\log\phi\left(\frac{y_i-\eta_0(\bx_i)}{\sigma_0}\right)$ is Lipschitz continuous.

To prove continuity of the term $\frac{1}{n}\sum_{i=1}^n\log\phi\left(\frac{y_i-\eta(\bx_i)}{\sigma}\right)$, we first recall from (A7) that $\log\phi(x)=\log\phi(|x|)$ 
is Lipschitz continuous in $x$.
Hence, if we can show that for each $i=1,\ldots,n$, $\frac{y_i-\eta(\bx_i)}{\sigma}$ is continuous in $(\by_n,\theta)$, then this would prove continuity of
$\frac{1}{n}\sum_{i=1}^n\log\phi\left(\frac{y_i-\eta(\bx_i)}{\sigma}\right)$ since sum and composition of continuous functions are continuous. 
Now, $|(y_{1i}-\eta_1(\bx_i))-(y_{2i}-\eta_2(\bx_i))|\leq|y_{1i}-y_{2i}|+|\eta_1(\bx_i)-\eta(\bx_i)|\leq\|\by_{1n}-\by_{2n}\|+\|\eta_1-\eta_2\|$, showing continuity
of $y_i-\eta(\bx_i)$. Division of this term by $\sigma~(>0)$, preserves continuity.

Hence, $f(\by_n,\theta)$ is continuous with respect to $(\by_n,\theta)$, so that (S1) holds in our case.

\subsection{Verification of (S2) and proof of Lemma \ref{lemma:nongp_lemma1}}
\label{subsec:nongp_S2}
It follows from (\ref{eq:nongp_like1}) and (\ref{eq:nongp_true_like1}), that
\begin{equation}
E_{\theta_0}\left[\frac{1}{n}\log\frac{f_{\theta_0}(\by_n)}{f_{\theta}(\by_n)}\right]=
\log\left(\frac{\sigma}{\sigma_0}\right)+\frac{1}{n}\sum_{i=1}^nE_{\theta_0}\left[\log\phi\left(\frac{y_i-\eta_0(\bx_i)}{\sigma_0}\right)\right]
-\frac{1}{n}\sum_{i=1}^nE_{\theta_0}\left[\log\phi\left(\frac{y_i-\eta(\bx_i)}{\sigma}\right)\right].
\label{eq:nongp_logratio1}
\end{equation}
Now $E_{\theta_0}\left[\log\phi\left(\frac{y_i-\eta_0(\bx_i)}{\sigma_0}\right)\right]=\int_{-\infty}^{\infty}\left[\log\phi(z)\right]\phi(z)dz=c$ (say),
so that for any $n\geq 1$,
\begin{equation}
\frac{1}{n}\sum_{i=1}^nE_{\theta_0}\left[\log\phi\left(\frac{y_i-\eta_0(\bx_i)}{\sigma_0}\right)\right]=c.
\label{eq:term1}
\end{equation}
Now for any $\bx\in\mathcal X$, let
\begin{equation}
g_{\eta,\sigma}(\bx)=E_{\theta_0}\left[\log\phi\left(\frac{y-\eta(\bx)}{\sigma}\right)\right]=\int_{-\infty}^{\infty}\log\phi\left(\frac{\sigma_0z+\eta_0(\bx)-\eta(\bx)}{\sigma}\right)\phi(z)dz.
\label{eq:g1}
\end{equation}
Let us first investigate continuity of $g_{\eta,\sigma}(\bx)$ with respect to $\bx$. To this end, observe that for $\bx_1,\bx_2\in\mathcal X$, the following hold 
thanks to Lipschitz continuity of $\log\phi$:
\begin{align}
&\left|g_{\eta,\sigma}(\bx_1)-g_{\eta,\sigma}(\bx_2)\right|\notag\\
&\ \ \leq\int_{-\infty}^{\infty}
\left|\log\phi\left(\frac{\sigma_0z+\eta_0(\bx_1)-\eta(\bx_1)}{\sigma}\right)-\log\phi\left(\frac{\sigma_0z+\eta_0(\bx_2)-\eta(\bx_2)}{\sigma}\right)\right|\phi(z)dz\notag\\
&\ \ =\frac{L}{\sigma}\int_{-\infty}^{\infty}\left|(\eta_0(\bx_1)-\eta_0(\bx_2))-(\eta(\bx_1)-\eta(\bx_2))\right|\phi(z)dz\notag\\
&\ \ \leq \frac{L}{\sigma}\left(\left|\eta_0(\bx_1)-\eta_0(\bx_2)\right|+\left|\eta(\bx_1)-\eta(\bx_2)\right|\right).\label{eq:g2}
\end{align}
In our model, $\eta(\bx)$ is continuous in $\bx$, but $\eta_0(\bx)$ need not be so. If $\eta_0(\bx)$ is allowed to be continuous, then by (\ref{eq:g2}), $g_{\eta,\sigma}(\bx)$
is continuous in $\bx$. If $\eta_0(\bx)$ has at most countably many discontinuities, then $g_{\eta,\sigma}(\bx)$ is continuous everywhere on $\mathcal X$
except perhaps at a countable number of points. In both the cases, $g_{\eta,\sigma}(\bx)$ is Riemann integrable when the covariates are considered deterministic. 
In that case, 
\begin{equation}
\frac{1}{n}\sum_{i=1}^nE_{\theta_0}\left[\log\phi\left(\frac{y_i-\eta(\bx_i)}{\sigma}\right)\right]=\frac{1}{n}\sum_{i=1}^ng_{\eta,\sigma}(\bx_i)\rightarrow
\int_{\mathcal X}g_{\eta,\sigma}(\bx)d\bx,~\mbox{as}~n\rightarrow\infty.
\label{eq:g3}
\end{equation}
If $\left\{\bx_i:i=1,2,\ldots\right\}$ is considered to be an $iid$ realization from $Q$, then by the ergodic theorem 
\begin{equation}
\frac{1}{n}\sum_{i=1}^nE_{\theta_0}\left[\log\phi\left(\frac{y_i-\eta(\bx_i)}{\sigma}\right)\right]=\frac{1}{n}\sum_{i=1}^ng_{\eta,\sigma}(\bx_i)\rightarrow
\int_{\mathcal X}g_{\eta,\sigma}(\bx)dQ,~\mbox{as}~n\rightarrow\infty.
\label{eq:g4}
\end{equation}
We denote both $\int_{\mathcal X}g_{\eta,\sigma}(\bx)d\bx$ and $\int_{\mathcal X}g_{\eta,\sigma}(\bx)dQ$ by $E_{\bX}\left[g_{\eta,\sigma}(\bX)\right]$.
Note that both the integrals exist thanks to continuity of $g_{\eta,\sigma}(\bx)$ and compactness of $\mathcal X$.
Combining (\ref{eq:term1}), (\ref{eq:g3}) and (\ref{eq:g4}) we obtain
\begin{equation}
E_{\theta_0}\left[\frac{1}{n}\log\frac{f_{\theta_0}(\by_n)}{f_{\theta}(\by_n)}\right]\rightarrow h(\theta),
\label{eq:g5}
\end{equation}
where $h(\theta)$ is given by (\ref{eq:nongp_h}).
In other words, (S2) holds.

\subsection{Verification of (S3) and proof of Theorem \ref{theorem:nongp_theorem1}}
\label{subsec:nongp_S3}
For any $\delta>0$, and for any $\theta\in\Theta$,
\begin{align}
&P\left(\left|\frac{1}{n}\log R_n(\theta)+h(\theta)\right|>\delta\right)\notag\\
&=P\left(\left|\frac{1}{n}\sum_{i=1}^n\log\phi\left(\frac{y_i-\eta(\bx_i)}{\sigma}\right)-\frac{1}{n}\sum_{i=1}^n\log\phi\left(\frac{y_i-\eta_0(\bx_i)}{\sigma_0}\right)
+c-E_{\bX}\left[g_{\eta,\sigma}(\bX)\right]\right|>\delta\right)\notag\\
&\leq P\left(\left|\frac{1}{n}\sum_{i=1}^n\log\phi\left(\frac{y_i-\eta(\bx_i)}{\sigma}\right)-E_{\bX}\left[g_{\eta,\sigma}(\bX)\right]\right|>\frac{\delta}{2}\right)
\label{eq:s3_1}\\
&\qquad + P\left(\left|\frac{1}{n}\sum_{i=1}^n\log\phi\left(\frac{y_i-\eta_0(\bx_i)}{\sigma_0}\right)-c\right|>\frac{\delta}{2}\right).
\label{eq:s3_2}
\end{align}
Let us focus attention on the probability given by (\ref{eq:s3_1}).
\begin{align}
&P\left(\left|\frac{1}{n}\sum_{i=1}^n\log\phi\left(\frac{y_i-\eta(\bx_i)}{\sigma}\right)-E_{\bX}\left[g_{\eta,\sigma}(\bX)\right]\right|>\frac{\delta}{2}\right)\notag\\
&=P\left(\left|\frac{1}{n}\sum_{i=1}^n\left[\log\phi\left(\frac{y_i-\eta(\bx_i)}{\sigma}\right)-g_{\eta,\sigma}(\bx_i)\right]
+\left[\frac{1}{n}\sum_{i=1}^ng_{\eta,\sigma}(\bx_i)-E_{\bX}\left[g_{\eta,\sigma}(\bX)\right]\right]\right|>\frac{\delta}{2}\right)\notag\\
&\leq P\left(\left|\frac{1}{n}\sum_{i=1}^n\left[\log\phi\left(\frac{y_i-\eta(\bx_i)}{\sigma}\right)-g_{\eta,\sigma}(\bx_i)\right]\right|>\frac{\delta}{4}\right)
\label{eq:s3_4}\\
&\qquad\qquad +P\left(\left|\frac{1}{n}\sum_{i=1}^ng_{\eta,\sigma}(\bx_i)-E_{\bX}\left[g_{\eta,\sigma}(\bX)\right]\right|>\frac{\delta}{4}\right).
\label{eq:s3_5}
\end{align}
Let us first deal with the probability given by (\ref{eq:s3_4}), with
$U_i=\log\phi\left(\frac{y_i-\eta(\bx_i)}{\sigma}\right)-g_{\eta,\sigma}(\bx_i)$.
Due to (A8), we apply Bernstein's inequality 
to obtain
\begin{align}
P\left(\left|\frac{1}{n}\sum_{i=1}^nU_i\right|>\frac{\delta}{4}\right)
&\leq 2\max\left\{P\left(\frac{1}{n}\sum_{i=1}^nU_i>\frac{\delta}{4}\right),P\left(\frac{1}{n}\sum_{i=1}^nU_i< -\frac{\delta}{4}\right)\right\}\notag\\
&\leq 2\exp\left(-\frac{n}{2}\min\left\{\frac{\delta^2}{16{s^2_{\eta,\sigma}}},\frac{\delta}{4s_{\eta,\sigma}}\right\}\right).
\label{eq:a7_5}
\end{align}

Now note that the probability given by (\ref{eq:s3_5}) is the probability of a deterministic quantity with respect to $\by_n$ and due to (\ref{eq:g3}) and (\ref{eq:g4}), 
is identically zero for large enough $n$. In the case of random covariates, using (A9) we obtain 
\begin{align}
\left|g_{\eta,\sigma}(\bx)\right|&\leq\int_{-\infty}^{\infty}\left|\log\phi\left(\frac{\sigma_0}{\sigma}z\right)\right|\phi(z)dz+\frac{L\|\eta-\eta_0\|}{\sigma}\notag\\
&\leq \frac{c_3+L\|\eta-\eta_0\|}{\sigma}=\tilde c_{\eta,\sigma}~\mbox{(say)}.
\label{eq:s7_6}
\end{align}
$g_{\eta,\sigma}(\bx_i)$ are independent, and satisfy (\ref{eq:s7_6}). Hence, Hoeffding's inequality 
yields
\begin{align}
&P\left(\left|\frac{1}{n}\sum_{i=1}^ng_{\eta,\sigma}(\bx_i)-E_{\bX}\left[g_{\eta,\sigma}(\bX)\right]\right|>\frac{\delta}{4}\right)\notag\\
&\qquad\leq \exp\left\{-\frac{n^2\delta^2}{144n\tilde c^2_{\eta,\sigma}}\right\}=\exp\left\{-\frac{n\delta^2}{144\tilde c^2_{\eta,\sigma}}\right\}.
\label{eq:a7_6}
\end{align}

The probability given by (\ref{eq:s3_2}) can be bounded in the same way as (\ref{eq:a7_5}). Indeed, we have
\begin{equation}
P\left(\left|\frac{1}{n}\sum_{i=1}^n\log\phi\left(\frac{y_i-\eta_0(\bx_i)}{\sigma_0}\right)-c\right|>\frac{\delta}{2}\right)
\leq 2\exp\left(-\frac{n}{2}\min\left\{\frac{\delta^2}{16{s^2_{\eta_0,\sigma_0}}},\frac{\delta}{4s_{\eta_0,\sigma_0}}\right\}\right).
\label{eq:a7_7}
\end{equation}

Combining the above results, it is seen that for any $\delta>0$, and for each $\theta\in\Theta$, there exists $a_{\theta}>0$, depending on $\theta$ such that
$P\left(\left|\frac{1}{n}\log R_n(\theta)+h(\theta)\right|>\delta\right)\leq 5\exp\left\{-na_{\theta}\right\}$, which is summable. Hence, by the Borel-Cantelli lemma,
$\frac{1}{n}\log R_n(\theta)\rightarrow -h(\theta)$, almost surely, as $n\rightarrow\infty$, for all $\theta\in\Theta$.
Thus, (S3) holds.

\subsection{Verification of (S4)}
\label{subsec:nongp_S4}

Using (\ref{eq:s7_6}) it is easily seen that 
\begin{equation}
h(\theta)\leq\left|\log\left(\frac{\sigma}{\sigma_0}\right)\right|+|c|+\int_{-\infty}^{\infty}\left|\log\phi\left(\frac{\sigma_0}{\sigma}z\right)\right|\phi(z)dz
+\frac{\|\eta_0\|+\|\eta\|}{\sigma}.
\label{eq:s4_1}
\end{equation}
Since almost surely with respect to the prior $\pi_{\sigma}$, $0<\sigma<\infty$, and $\|\eta\|<\infty$ almost surely with respect to the prior of $\eta$, and since $\|\eta_0\|<\infty$, 
it follows from (\ref{eq:s4_1}), 
that $\pi\left(h(\theta)=\infty\right)=0$, showing that (S4) holds.

\subsection{Verification of (S5)}
\label{subsec:nongp_S5}

\subsubsection{Verification of (S5) (1)}
\label{subsubsec:nongp_S5_1}


Recall from (\ref{eq:G}) that
\begin{align}
	\mathcal G_n&=\left\{\left(\eta,\sigma\right):\|\eta\|\leq\exp(\left(\beta n\right)^{1/4}),\exp(-\left(\beta n\right)^{1/4})\leq\sigma\leq\exp(\left(\beta n\right)^{1/4}),
	\|\eta'_j\|\leq\exp(\left(\beta n\right)^{1/4});j=1,\ldots,d\right\}.\notag
\end{align}
Then $\mathcal G_n\rightarrow\Theta$, as $n\rightarrow\infty$.
Now note that
\begin{align}
	&\pi(\mathcal G_n)=\pi\left(\|\eta\|\leq\exp(\left(\beta n\right)^{1/4}),\exp(-\left(\beta n\right)^{1/4})\leq\sigma\leq\exp(\left(\beta n\right)^{1/4})\right)\notag\\
	&\ \ -\pi\left(\left\{\|\eta'_j\|\leq\exp(\left(\beta n\right)^{1/4});j=1,\ldots,d\right\}^c\right)\notag\\
	&=\pi\left(\|\eta\|\leq\exp(\left(\beta n\right)^{1/4}),\exp(-\left(\beta n\right)^{1/4})\leq\sigma\leq\exp(\left(\beta n\right)^{1/4})\right)\notag\\
	&\ \ -\pi\left(\bigcup_{j=1}^d\left\{\|\eta'_j\|>\exp(\left(\beta n\right)^{1/4})\right\}\right)\notag\\
	& \geq 1-\pi\left(\|\eta\|>\exp(\left(\beta n\right)^{1/4})\right)-\pi\left(\left\{\exp(-\left(\beta n\right)^{1/4})
	\leq\sigma\leq\exp(\left(\beta n\right)^{1/4})\right\}^c\right)\notag\\
	&\ \ \  \ -\sum_{j=1}^d\pi\left(\|\eta'_j\|>\exp(\left(\beta n\right)^{1/4})\right)\notag\\
&\geq 1-(c_{\eta}+c_{\sigma}+\sum_{j=1}^dc_{\eta^\prime_j})\exp(-\beta n),
\label{eq:nongp_s5_1}
\end{align}
by (A5) and (A6). 
In other words, (S5) (1) holds.

\subsubsection{Verification of (S5) (2)}
\label{subsec:nongp_S5_2}
We now show that (S5) (2), namely, convergence in (S3) is uniform in $\theta$ over $\mathcal G_n\setminus I$ holds. 
In our case, by (S4), $h(\theta)<\infty$ with probability one, so that $I=\emptyset$ and $\mathcal G_n\setminus I=\mathcal G_n$.
Since we have already proved in the context of (S3) that $\underset{n\rightarrow\infty}{\lim}~\frac{1}{n}\log R_n(\theta)=-h(\theta)$, almost surely, for all $\theta\in\Theta$,
(S5) (2) will be verified if we can further prove that $\mathcal G_n$ is compact for each $n\geq 1$ and if $\frac{1}{n}\log R_n(\theta)+h(\theta)$ is Lipschitz in $\theta\in\mathcal G$,
for any $\mathcal G\in\left\{\mathcal G_n:n=1,2,\ldots\right\}$.

Compactness of $\mathcal G_n$, for all $n\geq 1$, follows as before.
Hence uniform convergence as required will be proven if we can show that $\frac{1}{n}\log R_n(\theta)+h(\theta)$ is 
stochastically equicontinuous almost surely in $\theta\in\mathcal G$ for any 
$\mathcal G\in\left\{\mathcal G_n:n=1,2,\ldots\right\}$ and $\frac{1}{n}\log R_n(\theta)+h(\theta)\rightarrow 0$, almost surely, for all $\theta\in\mathcal G$
As before, we rely on Lipschitz continuity.
To see that $\frac{1}{n}\log R_n(\theta)+h(\theta)$ is Lipschitz in $\theta\in\mathcal G$, first observe that it follows from the arguments in Section \ref{subsec:nongp_S1} that
$\frac{1}{n}\log R_n(\theta)$ is Lipschitz in $\eta$, when the data are held constant. Moreover, the derivative with respect to $\sigma$ is bounded since $\log\phi$ is Lipschitz
and since $\sigma$ in bounded in $\mathcal G$. In other words, $\frac{1}{n}\log R_n(\theta)$ is almost surely Lipschitz in $\theta$. Thus, if we can show that $h(\theta)$ is also Lipschitz in $\theta$,
then this would prove that $\frac{1}{n}\log R_n(\theta)+h(\theta)$ is almost surely Lipschitz in $\theta$. For our purpose, it is sufficient to show that $E_{\bX}\left[g_{\eta,\sigma}(\bX)\right]$
is Lipschitz in $(\eta,\sigma)$. Since for any $\eta_1,\eta_2,\sigma\in\Theta$, 
\begin{align}
&\left|E_{\bX}\left[g_{\eta_1,\sigma}(\bX)\right]-E_{\bX}\left[g_{\eta_2,\sigma}(\bX)\right]\right|\notag\\
&\qquad\leq E_{\bX}\left|g_{\eta_1,\sigma}(\bX)-g_{\eta_2,\sigma}(\bX)\right|\notag\\
&\qquad =E_{\bX}\left[\int_{-\infty}^{\infty}\left|\log\phi\left(\frac{\sigma_0 z+\eta_0(\bX)-\eta_1(\bX)}{\sigma}\right)
-\log\phi\left(\frac{\sigma_0 z+\eta_0(\bX)-\eta_2(\bX)}{\sigma}\right)\right|\phi(z)dz\right]\notag\\
&\qquad\leq\frac{L}{\sigma}E_{\bX}\left[\int_{-\infty}^{\infty}\left|\eta_1(\bX)-\eta_2(\bX)\right|\phi(z)dz\right]\notag\\
&\qquad\leq\frac{L_2}{\sigma}\|\eta_1-\eta_2\|,
\label{eq:lip1}
\end{align}
for some $L_2>0$, $E_{\bX}\left[g_{\eta,\sigma}(\bX)\right]$ is Lipschitz in $\eta$.
Now recall that under the assumption (A9), $\int_{-\infty}^{\infty}|z|\phi(z)dz<\infty$.
With this, we now show that $E_{\bX}\left[g_{\eta,\sigma}(\bX)\right]$ has bounded first derivative with respect to $\sigma$ in the interior of $\mathcal G$.
Observe that
\begin{align}
&|r|^{-1}\left|g_{\eta,\sigma+r}(\bx)-g_{\eta,\sigma}(\bx)\right|\notag\\
&\qquad\leq |r|^{-1}\left[\int_{-\infty}^{\infty}\left|\log\phi\left(\frac{\sigma_0 z+\eta_0(\bx)-\eta(\bx)}{\sigma+r}\right)
-\log\phi\left(\frac{\sigma_0 z+\eta_0(\bx)-\eta(\bx)}{\sigma}\right)\right|\phi(z)dz\right]\notag\\
&\qquad\leq L \left[\int_{-\infty}^{\infty}\left(\frac{\sigma_0|z|+\|\eta-\eta_0\|}{\sigma (\sigma+r)}\right)\phi(z)dz\right]~\left(\mbox{since}~\log\phi~\mbox{is Lipschitz}\right)\notag\\
&\qquad\leq \frac{L}{\sigma (\sigma+r)}\left(\sigma_0\int_{-\infty}^{\infty}|z|\phi(z)dz+\|\eta-\eta_0\|\right).
\label{eq:lip2}
\end{align}
By (A9), $\int_{-\infty}^{\infty}|z|\phi(z)dz<\infty$, and $\sigma$, $\sigma+r$ (both in the interior of $\mathcal G$) are both upper and lower bounded in $\mathcal G$, the lower bound
being strictly positive. 
Hence, (\ref{eq:lip2}) is integrable with respect to (the distribution) of $\bX$, since $\mathcal X$ is compact. 
Hence, by the dominated convergence theorem, differentiation with respect to $\sigma$ can be performed inside the double integral associated with $E_{\bX}\left[g_{\eta,\sigma}(\bX)\right]$.
Since $\log\phi$ has bounded first derivative as it is Lipschitz and since $\sigma$ is lower bounded by a positive quantity in $\mathcal G$, it 
follows that $E_{\bX}\left[g_{\eta,\sigma}(\bX)\right]$ has bounded 
first derivative with respect to $\sigma$. Combined with the result that $E_{\bX}\left[g_{\eta,\sigma}(\bX)\right]$ is Lipschitz in $\eta$, this yields that 
$E_{\bX}\left[g_{\eta,\sigma}(\bX)\right]$ is Lipschitz in $(\eta,\sigma)$. In conjunction with the result that $\frac{1}{n}\log R_n(\theta)$ is almost surely Lipschitz in $\theta$, it holds that
$\frac{1}{n}\log R_n(\theta)+h(\theta)$ is almost surely Lipschitz in $\theta\in\mathcal G$.
In other words, (S5) (2) stands verified.

\subsubsection{Verification of (S5) (3)}
\label{subsec:nongp_S5_3}
To verify (S5) (3), note that continuity of $h(\theta)$, compactness of $\mathcal G_n$, along with its non-decreasing nature with respect to $n$ implies that 
$h\left(\mathcal G_n\right)\rightarrow h(\Theta)$, as $n\rightarrow\infty$.

\subsection{Verification of (S6) and proof of Theorem \ref{theorem:nongp_theorem3}}
\label{subsec:nongp_S6}

Let $\kappa_1=\kappa-h(\Theta)$. 
%
Then 
it follows from (\ref{eq:s3_2}), (\ref{eq:s3_4}), (\ref{eq:s3_5}), (\ref{eq:a7_5}), (\ref{eq:a7_6}) and (\ref{eq:a7_7}),
that
\begin{align}
&\int_{\mathcal S^c}P\left(\left|\frac{1}{n}\log R_n(\theta)+h(\theta)\right|>\kappa_1\right)d\pi(\theta)\notag\\
&\leq \int_{\mathcal S^c}P\left(\left|\frac{1}{n}\sum_{i=1}^n\left[\log\phi\left(\frac{y_i-\eta(\bx_i)}{\sigma}\right)-g_{\eta,\sigma}(\bx_i)\right]\right|>\frac{\kappa_1}{4}\right)d\pi(\theta)\notag\\
&\qquad+\int_{\mathcal S^c}P\left(\left|\frac{1}{n}\sum_{i=1}^ng_{\eta,\sigma}(\bx_i)-E_{\bX}\left[g_{\eta,\sigma}(\bX)\right]\right|>\frac{\kappa_1}{4}\right)d\pi(\theta)\notag\\
&\qquad + \int_{\mathcal S^c}P\left(\left|\frac{1}{n}\sum_{i=1}^n\log\phi\left(\frac{y_i-\eta_0(\bx_i)}{\sigma_0}\right)-c\right|>\frac{\kappa_1}{2}\right)d\pi(\theta)\notag\\
&\leq \int_{\mathcal S^c}2\exp\left(-\frac{n}{2}\min\left\{\frac{\kappa^2_1}{16s^2_{\eta,\sigma}},\frac{\kappa_1}{4s_{\eta,\sigma}}\right\}\right)d\pi(\theta)
+\int_{\mathcal S^c}\exp\left(-\frac{Cn\kappa^2_1}{\tilde c^2_{\eta,\sigma}}\right)d\pi(\theta)
\label{eq:s6_1_supp}\\
&\qquad+\int_{\mathcal S^c}2\exp\left(-\frac{n}{2}\min\left\{\frac{\kappa^2_1}{16s^2_{\eta_0,\sigma_0}},\frac{\kappa_1}{4s_{\eta_0,\sigma_0}}\right\}\right)d\pi(\theta),
\label{eq:s6_2_supp}
\end{align}
for some relevant positive constant $C$.

Now, in the same way as (\ref{eq:summable3}) we obtain
\begin{align}
&\int_{\mathcal S^c}2\exp\left(-\frac{n}{2}\min\left\{\frac{\kappa^2_1}{16s^2_{\eta,\sigma}},\frac{\kappa_1}{4s_{\eta,\sigma}}\right\}\right)d\pi(\theta)\notag\\
&\leq \int_{\mathcal G_n}2\exp\left(-\frac{n}{2}\frac{\kappa^2_1}{16s^2_{\eta,\sigma}}\right)d\pi(\theta)
+ \int_{\mathcal G_n}2\exp\left(-\frac{n}{2}\frac{\kappa_1}{4s_{\eta,\sigma}}\right)d\pi(\theta)+2\pi(\mathcal G^c_n)\notag\\
	&\leq C_1\exp\left\{-\left(C_2\sqrt{\kappa_1}n^{1/4}-5\left(\beta n\right)^{1/4}-2n^q\log c_5\right)\right\}\notag\\
	&\qquad+\tilde C_1\exp\left\{-\left(\tilde C_2\sqrt{\kappa_1}n^{1/4}-\frac{9}{2}\left(\beta n\right)^{1/4}-2n^q\log c_5\right)\right\}+2\pi\left(\mathcal G^c_n\right),
\label{eq:s6_4_supp}
\end{align}
for relevant positive constants $C_1,C_2,\tilde C_1,\tilde C_2,c_5$. In the same way,
\begin{equation}
\int_{\mathcal S^c}\exp\left(-\frac{Cn\kappa^2_1}{\tilde c^2_{\eta,\sigma}}\right)d\pi(\theta)
\leq C_{11}\exp\left\{-\left(C_{21}\sqrt{\kappa_1}n^{1/4}-5\left(\beta n\right)^{1/4}-2n^q\log c_{51}\right)\right\} +\pi\left(\mathcal G^c_n\right),
	\label{eq:s6_5_supp}
\end{equation}
for relevant positive constants $C_{11},C_{21},c_{51}$.

Since (\ref{eq:s6_4_supp}) and (\ref{eq:s6_5_supp}) are summable, and since (\ref{eq:s6_2_supp}) is summable (as the integrand is independent of parameters),
it follows from (\ref{eq:s6_1_supp}) and (\ref{eq:s6_2_supp}) that
\begin{equation*}
\int_{\mathcal S^c}P\left(\left|\frac{1}{n}\log R_n(\theta)+h(\theta)\right|>\kappa_1\right)d\pi(\theta)<\infty,
\end{equation*}
showing that (S6) holds.

\subsection{Verification of (S7)}
\label{subsec:nongp_S7}
For any set $A$ such that $\pi(A)>0$, $\mathcal G_n\cap A\uparrow A$. It follows from this and continuity of $h$ that $h\left(\mathcal G_n\cap A\right)\downarrow h\left(A\right)$ as
$n\rightarrow\infty$, so that (S7) holds.

\bibliography{irmcmc}

\begin{thebibliography}{}

\bibitem[Adler(1981)Adler]{Adler81}
Adler, R.~J. (1981).
\newblock {\em {T}he {G}eometry of {R}andom {F}ields\/}.
\newblock John Wiley \& Sons Ltd., New York.

\bibitem[Adler and Taylor(2007)Adler and Taylor]{Adler07}
Adler, R.~J. and Taylor, J.~E. (2007).
\newblock {\em {R}andom {F}ields and {G}eometry\/}.
\newblock Springer, New York.

\bibitem[Banerjee {\em et~al.}(2014)Banerjee, Carlin, and Gelfand]{Banerjee14}
Banerjee, S., Carlin, B.~P., and Gelfand, A.~E. (2014).
\newblock {\em {H}ierarchical {M}odeling and {A}nalysis for {S}patial
  {D}ata\/}.
\newblock Chapman \& Hall/CRC, USA.

\bibitem[Bennett(1962)Bennett]{Bennett62}
Bennett, G. (1962).
\newblock {P}robability {I}nequalities for the {S}ums of {I}ndependent {R}andom
  {V}ariables.
\newblock {\em Journal of the American Statistical Association\/}, {\bf 57},
  33--45.

\bibitem[Billingsley(2013)Billingsley]{Billingsley13}
Billingsley, P. (2013).
\newblock {\em {C}onvergence of {P}robability {M}easures\/}.
\newblock John Wiley \& Sons, New Jersey, USA.

\bibitem[Choi(2009)Choi]{Choi09}
Choi, T. (2009).
\newblock {A}symptotic {P}roperties of {P}osterior {D}istributions in
  {N}onparametric {R}egression with {N}on-{G}aussian {E}rrors.
\newblock {\em Annals of the Institute of Statistical Mathematics\/}, {\bf 61},
  835--859.

\bibitem[Choi and Schervish(2007)Choi and Schervish]{Choi07}
Choi, T. and Schervish, M.~J. (2007).
\newblock {O}n {P}osterior {C}onsistency in {N}onparametric {R}egression
  {P}roblems.
\newblock {\em Journal of Multivariate Analysis\/}, {\bf 98}, 1969--1987.

\bibitem[Cramer and Leadbetter(1967)Cramer and Leadbetter]{Cramer67}
Cramer, H. and Leadbetter, M.~R. (1967).
\newblock {\em {S}tationary and {R}elated {S}tochastic {P}rocesses\/}.
\newblock Wiley, New York.

\bibitem[Cressie(1993)Cressie]{Cressie93}
Cressie, N. A.~C. (1993).
\newblock {\em {S}tatistics for {S}patial {D}ata\/}.
\newblock Wiley, New York.

\bibitem[Dey {\em et~al.}(2012)Dey, M{\"u}ller, and Sinha]{Dey12}
Dey, D.~K., M{\"u}ller, P., and Sinha, D. (2012).
\newblock {\em {P}ractical {N}onparametric and {S}emiparametric {B}ayesian
  {S}tatistics\/}.
\newblock Springer, New York, USA.

\bibitem[Efromovich(2008)Efromovich]{Efro08}
Efromovich, S. (2008).
\newblock {\em {N}onparametric {C}urve {E}stimation: {M}ethods, {T}heory, and
  {A}pplications\/}.
\newblock Springer, New York, USA.

\bibitem[Eubank(1999)Eubank]{Eubank99}
Eubank, R. (1999).
\newblock {\em {N}onparametric {R}egression and {S}pline {S}moothing\/}.
\newblock CRC Press, New York, USA.

\bibitem[Ghosal and {van derVaart}(2017)Ghosal and {van derVaart}]{Ghosal17}
Ghosal, A. and {van derVaart}, A. (2017).
\newblock {\em {F}undamentals of {N}onparametric {B}ayesian {I}nference\/}.
\newblock Cambridge University Press, Cambridge, UK.

\bibitem[Ghosh and Ramamoorthi(2003)Ghosh and Ramamoorthi]{Ghosh03}
Ghosh, J.~K. and Ramamoorthi, R.~V. (2003).
\newblock {\em {B}ayesian {N}onparametrics\/}.
\newblock Springer, New York, USA.

\bibitem[Green and Silverman(1993)Green and Silverman]{Green93}
Green, P.~J. and Silverman, B.~W. (1993).
\newblock {\em {N}onparametric {R}egression and {G}eneralized {L}inear
  {M}odels: {A} {R}oughness {P}enalty {A}pproach\/}.
\newblock CRC Press, New York, USA.

\bibitem[H{\"a}rdle(1990)H{\"a}rdle]{Hardle90}
H{\"a}rdle, W.~K. (1990).
\newblock {\em {A}pplied {N}onparametric {R}egression\/}.
\newblock Cambridge University Press, Cambridge, UK.

\bibitem[H{\"a}rdle {\em et~al.}(2012)H{\"a}rdle, M{\"u}ller, Sperlich, and
  Werwatz]{Hardle12}
H{\"a}rdle, W.~K., M{\"u}ller, M., Sperlich, S., and Werwatz, A. (2012).
\newblock {\em {N}onparametric and {S}emiparametric {M}odels\/}.
\newblock Springer, New York, USA.

\bibitem[Hjort {\em et~al.}(2010)Hjort, Holmes, M{\"u}ller, and
  Walker]{Hjort10}
Hjort, N.~L., Holmes, C., M{\"u}ller, P., and Walker, S.~G. (2010).
\newblock {\em {B}ayesian {N}onparametrics\/}.
\newblock Cambridge University Press, Cambridge, UK.

\bibitem[Hoeffding(1963)Hoeffding]{Hoeffding63}
Hoeffding, W. (1963).
\newblock {P}robability {I}nequalities for {S}ums of {B}ounded {R}andom
  {V}ariables.
\newblock {\em Journal of the American Statistical Association\/}, {\bf 58},
  13--30.

\bibitem[Knapik and Salomond(2018)Knapik and Salomond]{Knapik18}
Knapik, B. and Salomond, J.~B. (2018).
\newblock {A} {G}eneral {A}pproach to {P}osterior {C}ontraction in
  {N}onparametric {I}nverse {P}roblems.
\newblock {\em Bernoulli\/}, {\bf 24}, 2091--2121.

\bibitem[Knapik {\em et~al.}(2011)Knapik, {van der Vaart}, and {van
  Zanten}]{Knapik11}
Knapik, B.~T., {van der Vaart}, A.~W., and {van Zanten}, J.~H. (2011).
\newblock {B}ayesian {I}nverse {P}roblems with {G}aussian {P}riors.
\newblock {\em The Annals of Statistics\/}, {\bf 39}, 2626--2657.

\bibitem[Lange(2010)Lange]{Lange10}
Lange, K. (2010).
\newblock {\em {N}umerical {A}nalysis for {S}tatisticians\/}.
\newblock New York, Springer.

\bibitem[Lenk(1988)Lenk]{Lenk88}
Lenk, P.~J. (1988).
\newblock {T}he {L}ogistic {N}ormal {D}istribution for {B}ayesian,
  {N}onparametric, {P}redictive {D}ensities.
\newblock {\em Journal of the American Statistical Association\/}, {\bf 83},
  509--516.

\bibitem[Lenk(1991)Lenk]{Lenk91}
Lenk, P.~J. (1991).
\newblock {T}owards a {P}racticable {B}ayesian {N}onparametric {D}ensity
  {E}stimator.
\newblock {\em Biometrika\/}, {\bf 78}, 531--543.

\bibitem[Lenk(2003)Lenk]{Lenk03}
Lenk, P.~J. (2003).
\newblock {B}ayesian {S}emiparametric {D}ensity {E}stimation and {M}odel
  {V}erification {U}sing a {L}ogistic-{G}aussian {P}rocess.
\newblock {\em Journal of Computational and Graphical Statistics\/}, {\bf 12},
  548--565.

\bibitem[Massart(2003)Massart]{Massart03}
Massart, P. (2003).
\newblock {C}oncentration {I}nequalities and {M}odel {S}election.
\newblock Volume 1896 of Lecture Notes in Mathematics. Springer-Verlag.
  Lectures given at the 33rd Probability Summer School in Saint-Flour.

\bibitem[M{\"u}ller {\em et~al.}(2015)M{\"u}ller, Quintana, Jara, and
  Hanson]{Muller15}
M{\"u}ller, P., Quintana, F.~A., Jara, A., and Hanson, T. (2015).
\newblock {\em {B}ayesian {N}onparametric {D}ata {A}nalysis\/}.
\newblock Springer, New York, USA.

\bibitem[Newey(1991)Newey]{Newey91}
Newey, W.~K. (1991).
\newblock {U}niform {C}onvergence in {P}robability and {S}tochastic
  {E}quicontinuity.
\newblock {\em Econometrica\/}, {\bf 59}, 1161--1167.

\bibitem[Rasmussen and Williams(2006)Rasmussen and Williams]{Ras06}
Rasmussen, C.~E. and Williams, C. K.~I. (2006).
\newblock {\em {G}aussian {P}rocesses for {M}achine {L}earning\/}.
\newblock MIT Press, Cambridge, MA.

\bibitem[Santner {\em et~al.}(2003)Santner, Williams, and Notz]{Santner03}
Santner, T.~J., Williams, B.~J., and Notz, W.~I. (2003).
\newblock {\em {T}he {D}esign and {A}nalysis of {C}omputer {E}xperiments\/}.
\newblock Springer Series in Statistics. Springer-Verlag, New York, Inc.

\bibitem[Schimek(2013)Schimek]{Schimek13}
Schimek, M.~J. (2013).
\newblock {\em {S}moothing and {R}egression: {A}pproaches, {C}omputation, and
  {A}pplication\/}.
\newblock John Wiley and Sons, New Jersey, USA.

\bibitem[Shalizi(2009)Shalizi]{Shalizi09}
Shalizi, C.~R. (2009).
\newblock {D}ynamics of {B}ayesian {U}pdating {W}ith {D}ependent {D}ata and
  {M}isspecified {M}odels.
\newblock {\em Electronic Journal of Statistics\/}, {\bf 3}, 1039--1074.

\bibitem[Takezawa(2006)Takezawa]{Tak06}
Takezawa, K. (2006).
\newblock {\em {I}ntroduction to {N}onparametric {R}egression\/}.
\newblock John Wiley \& Sons, New Jersey, USA.

\bibitem[Uspensky(1937)Uspensky]{Uspensky37}
Uspensky, J.~V. (1937).
\newblock {\em {I}ntroduction to {M}athematical {P}robability\/}.
\newblock McGraw-Hill, New York, USA.

\bibitem[{van der Vaart} and {van Zanten}(2008){van der Vaart} and {van
  Zanten}]{Vaart08}
{van der Vaart}, A.~W. and {van Zanten}, J.~H. (2008).
\newblock {R}ates of {C}ontraction of {P}osterior {D}istributions {B}ased on
  {G}aussian {P}rocess {P}riors.
\newblock {\em The Annals of Statistics\/}, {\bf 36}, 1435--1463.

\bibitem[{van der Vaart} and {van Zanten}(2009){van der Vaart} and {van
  Zanten}]{Vaart09}
{van der Vaart}, A.~W. and {van Zanten}, J.~H. (2009).
\newblock {A}daptive {B}ayesian {E}stimation {U}sing a {G}aussian {R}andom
  {F}ield with {I}nverse {G}amma {B}andwidth.
\newblock {\em The Annals of Statistics\/}, {\bf 37}, 2655--2675.

\bibitem[{van der Vaart} and {van Zanten}(2011){van der Vaart} and {van
  Zanten}]{Vaart11}
{van der Vaart}, A.~W. and {van Zanten}, J.~H. (2011).
\newblock {I}nformation {R}ates of {N}onparametric {G}aussian {P}rocess
  {M}ethods.
\newblock {\em Journal of Machine Learning Research\/}, {\bf 12}, 2095--2119.

\bibitem[Vollmer(2013)Vollmer]{Vollmer13}
Vollmer, S.~J. (2013).
\newblock {P}osterior {C}onsistency for {B}ayesian {I}nverse {P}roblems
  {T}hrough {S}tability and {R}egression {R}esults.
\newblock {\em Inverse Problems\/}, {\bf 29}, 125011.

\bibitem[Wu and Zhang(2006)Wu and Zhang]{Wu06}
Wu, H. and Zhang, J.-T. (2006).
\newblock {\em {N}onparametric {R}egression {M}ethods for {L}ongitudinal {D}ata
  {A}nalysis: {M}ixed-{E}ffects {M}odeling {A}pproaches\/}.
\newblock John Wiley and Sons, New Jersey, USA.

\bibitem[Yang {\em et~al.}(2018)Yang, Bhattacharya, and D.Pati]{yang17}
Yang, Y., Bhattacharya, A., and D.Pati (2018).
\newblock {F}requentist {C}overage and {S}up-norm {C}onvergence {R}ate in
  {G}aussian {P}rocess {R}egression.
\newblock ArXiv Preprint.

\end{thebibliography}


\begin{thebibliography}{}

\bibitem[Bennett(1962)Bennett]{Bennett62}
Bennett, G. (1962).
\newblock {P}robability {I}nequalities for the {S}ums of {I}ndependent {R}andom
  {V}ariables.
\newblock {\em Journal of the American Statistical Association\/}, {\bf 57},
  33--45.

\bibitem[Billingsley(2013)Billingsley]{Billingsley13}
Billingsley, P. (2013).
\newblock {\em {C}onvergence of {P}robability {M}easures\/}.
\newblock John Wiley \& Sons, New Jersey, USA.

\bibitem[Hoeffding(1963)Hoeffding]{Hoeffding63}
Hoeffding, W. (1963).
\newblock {P}robability {I}nequalities for {S}ums of {B}ounded {R}andom
  {V}ariables.
\newblock {\em Journal of the American Statistical Association\/}, {\bf 58},
  13--30.

\bibitem[Massart(2003)Massart]{Massart03}
Massart, P. (2003).
\newblock {C}oncentration {I}nequalities and {M}odel {S}election.
\newblock Volume 1896 of Lecture Notes in Mathematics. Springer-Verlag.
  Lectures given at the 33rd Probability Summer School in Saint-Flour.

\bibitem[Newey(1991)Newey]{Newey91}
Newey, W.~K. (1991).
\newblock {U}niform {C}onvergence in {P}robability and {S}tochastic
  {E}quicontinuity.
\newblock {\em Econometrica\/}, {\bf 59}, 1161--1167.

\bibitem[Shalizi(2009)Shalizi]{Shalizi09}
Shalizi, C.~R. (2009).
\newblock {D}ynamics of {B}ayesian {U}pdating {W}ith {D}ependent {D}ata and
  {M}isspecified {M}odels.
\newblock {\em Electronic Journal of Statistics\/}, {\bf 3}, 1039--1074.

\bibitem[Uspensky(1937)Uspensky]{Uspensky37}
Uspensky, J.~V. (1937).
\newblock {\em {I}ntroduction to {M}athematical {P}robability\/}.
\newblock McGraw-Hill, New York, USA.

\end{thebibliography}

\end{document}